\documentclass[12pt,reqno]{amsart}
\author{Kai (Steve) Fan}
\author{Paul Pollack}
\address{Department of Mathematics\\ University of Georgia\\ Athens, GA 30602}
\email{Steve.Fan@uga.edu}
\email{pollack@uga.edu}
\usepackage{mathtools}
\usepackage{empheq}
\usepackage{colonequals}

\newif\ifcoloringon
\coloringonfalse

\title{Counting primes with a given primitive root, uniformly}
\usepackage{amsmath,amssymb,amsthm,microtype,geometry,mathscinet,enumitem,url}
\usepackage{parskip}
\usepackage{stackengine}
\usepackage{colonequals}
\usepackage{xcolor}

\renewcommand\supset\supseteq
\renewcommand\epsilon\varepsilon
\renewcommand\phi\varphi
\usepackage[utf8]{inputenc}
\makeatletter
\NewCommandCopy\@@pmod\pmod
\DeclareRobustCommand{\pmod}{\@ifstar\@pmods\@@pmod}
\def\@pmods#1{\mkern4mu({\operator@font mod}\mkern 6mu#1)}
\makeatother

\numberwithin{equation}{section}

\DeclareMathAlphabet{\curly}{U}{rsfs}{m}{n}
\newtheorem{thm}{Theorem}[section]
\newtheorem{cor}[thm]{Corollary}
\newtheorem{prop}[thm]{Proposition}
\newtheorem{lem}[thm]{Lemma}

\newtheorem*{ALS}{Arithmetic large sieve}
\newtheorem*{GLS}{Larger sieve}

\theoremstyle{remark}
\newtheorem*{rmk}{Remark}

\usepackage{enumitem}
\coloringontrue

\begin{document}

\def\Cc{\curly{C}}
\def\Ll{\mathcal{L}}
\def\Ss{\mathcal{S}}
\def\Ee{\curly{E}}
\def\Oo{\mathcal{O}}
\def\Aa{\mathcal{A}}
\def\I{\mathcal{I}}
\def\J{\mathcal{J}}
\def\N{\mathbb{N}}
\def\Q{\mathbb{Q}}
\def\Qq{\mathcal{Q}}
\def\Rr{\mathcal{R}}
\def\Cl{\mathrm{Cl}}
\def\radi{\mathrm{rad}}
\def\PrinCl{\mathrm{PrinCl}}
\def\PreCl{\mathrm{PreCl}}
\newcommand\rad{\mathrm{radexp}}
\def\E{\mathbb{E}}
\def\Gg{\mathcal{G}}
\def\Z{\mathbb{Z}}
\def\F{\mathbb{F}}
\def\PP{\mathbb{P}}
\def\R{\mathbb{R}}
\def\C{\mathbb{C}}
\def\Gal{\mathrm{Gal}}
\def\Dd{\mathcal{D}}
\def\Pp{\mathcal{P}}
\def\Ss{\mathcal{S}}
\def\probp{\mathrm{prob}}
\def\omegab{\overline{\nu}}
\newcommand\Frob{\mathrm{Frob}}
\def\lcm{\mathop{\mathrm{lcm}}}
\renewcommand\subset\subseteq
\newcommand\Li{\mathrm{Li}}
\newcommand\altxrightarrow[2][0pt]{\mathrel{\ensurestackMath{\stackengine%
  {\dimexpr#1-7.5pt}{\xrightarrow{\phantom{#2}}}{\scriptstyle\!#2\,}%
  {O}{c}{F}{F}{S}}}}
\newcommand{\congto}{\altxrightarrow{\sim}}


\dedicatory{For Greg Martin on his retirement.}
\begin{abstract} 
The celebrated Artin conjecture on primitive roots asserts that given any integer $g$ which is neither $-1$ nor a perfect square, there is an explicit constant $A(g)>0$ such that the number $\Pi(x;g)$ of primes $p\le x$ for which $g$ is a primitive root is asymptotically $A(g)\pi(x)$ as $x\to\infty$, where $\pi(x)$ counts the number of primes not exceeding $x$. Artin's conjecture has remained unsolved since its formulation in 1927. Nevertheless, Hooley demonstrated in 1967 that Artin's conjecture is a consequence of the Generalized Riemann Hypothesis (GRH) for Dedekind zeta functions of certain cyclotomic-Kummer extensions over $\Q$. 
In this paper, we use  GRH  to establish a uniform version of the Artin--Hooley asymptotic formula. Specifically, we prove that $\Pi(x;g) \sim A(g) x/\log{x}$ whenever $\log{x}/\log\log{2|g|} \to \infty$, i.e., whenever $x$ tends to infinity faster than any power of $\log{(2|g|)}$. Under GRH, we also show that the least prime $p_g$ possessing $g$ as a primitive root satisfies the upper bound $p_g=O(\log^{19}(2|g|))$ uniformly for all non-square $g\ne-1$. We conclude with an application to the average value of $p_g$ and a discussion of an analogue concerning the least ``almost-primitive'' root.
\end{abstract}

\maketitle

\section{Introduction}
It is a classical result, due to Gauss, that the multiplicative group modulo a prime $p$ is always cyclic. That is, given any prime number $p$, there is an integer $g$ whose reduction mod $p$ generates the group $(\Z/p\Z)^{\times}$; following tradition, we call such an integer $g$ a \textsf{primitive root} modulo $p$. On the other hand, if we start with a given $g\in \Z$, there need not be any prime $p$ with $g$ a primitive root mod $p$. For instance, $g=4$ is not a primitive root modulo any prime, and the same holds for all even square values of $g$.

The distribution of primes $p$ possessing a prescribed integer $g$ as a primitive root is the subject of a celebrated 1927 conjecture of Emil Artin, formulated during a visit of Artin to Hasse (consult \cite[\S17.2]{AH08} for the history, and see \cite{moree12} for a comprehensive survey of related developments). For real $x>0$ and integers $g$, let
\[ \Pi(x;g) = \#\{\text{primes }p \le x: g\text{ is a primitive root mod $p$}\}.\]
Let \[ \Gg = \{g \in \Z: |g| > 1, g\text{ not a square}\}.\] \textsf{Artin's primitive root conjecture} predicts that for each $g \in \Gg$, \begin{equation}\label{eq:ACprecise} \Pi(x;g) \sim A(g) \pi(x), \quad\text{as}\quad x\to\infty, \end{equation} for an explicitly given constant $A(g) > 0$. 

The conjectured form of $A(g)$ depends on the arithmetic nature of $g$. For each $g \in \Gg$, let $g_1$ denote the unique squarefree integer with $g \in g_1 (\Q^{\times})^2$, and let $h$ be the largest positive integer for which $g \in (\Q^{\times})^h$. Since $g$ is not a square, $h$ is odd. Put
\begin{equation}\label{eq:A0eq} A_0(g) = \prod_{q \mid h}\bigg(1-\frac{1}{q-1}\bigg) \prod_{q \nmid h}\bigg(1-\frac{1}{q(q-1)}\bigg).\end{equation}
If $g_1 \equiv 1\pmod{4}$, put 
\begin{equation}\label{eq:A1eq} A_1(g) = 1-\mu(|g_1|) \prod_{\substack{q\mid h \\ q\mid g_1}}\frac{1}{q-2}\prod_{\substack{q\nmid h \\ q\mid g_1}}\frac{1}{q^2-q-1}; \end{equation}
otherwise, set $A_1(g) = 1$. Finally, put
\[ A(g) = A_0(g) A_1(g).\]
It is this value of $A(g)$ for which Artin predicts the asymptotic formula \eqref{eq:ACprecise}.\footnote{Artin's original 1927 formulation was missing the factor of $A_1(g)$. Artin realized the need for $A_1(g)$ after learning of computations carried out by the Lehmers. See Stevenhagen's discussion in \cite{stevenhagen03}.} 

Artin's conjecture remains unresolved. In fact, to this day there is not a single value of $g$ for which we can show even the weaker assertion that $\Pi(x;g) \to \infty$ as $x\to\infty$. (However, work of Heath-Brown \cite{HB86} implies this holds for at least one of $g=2$, $3$, or $5$.) The most important progress in this direction is a 1967 theorem of Hooley \cite{hooley67}, asserting that the full asymptotic relation \eqref{eq:ACprecise} follows from the Generalized Riemann Hypothesis (GRH).\footnote{Here and below, GRH means the Riemann Hypothesis for all Dedekind zeta functions of number fields.}

Hooley states and proves his asymptotic formula for \emph{fixed} $g \in \Gg$. Our main result makes the dependence on $g$ explicit.

\begin{thm}[assuming GRH]\label{thm:uniformasymptotic}The asymptotic formula $\Pi(x;g) \sim A(g) \pi(x)$ holds whenever $\log{x}/\log\log{2|g|} \to\infty$. More precisely, there is an absolute constant $x_0 > 0$ for which the following holds: If $g \in \Gg$ and $x \ge \max\{x_0, \log^3(2|g|) \}$, then
\begin{equation}\label{eq:uniformasymptotic} \Pi(x;g) = A(g)\pi(x) \left(1 + O\left(\frac{\log\log{x}}{\log{x}} + \frac{\log\log{2|g|}}{\log{x}}\right)\right).\end{equation}
\end{thm}

The proof of Theorem \ref{thm:uniformasymptotic}, presented in \S\ref{sec:uniformhooley}, broadly proceeds along the same course as Hooley's, but care and caution are required to ensure the final estimate is nontrivial in a wide range of $x$ and $g$. In particular, the fact that the positive constant $A(g)$ can be arbitrarily small causes substantial complications. 

Let $p_g$ denote the least prime $p$ possessing $g$ as a primitive root, where we set $p_g=\infty$ when no such $p$ exists. Theorem \ref{thm:uniformasymptotic} implies immediately that for all $g \in \Gg$,
\begin{equation}\label{eq:pgupper} p_g \ll \log^{B}({2|g|}),\end{equation}
for a certain absolute constant $B$. Indeed, suppose that $K$ is an admissible value of the implied constant in \eqref{eq:uniformasymptotic} and fix any constant $B> \max\{3,K\}$. If $x\ge \max\{x_0,\log^{B}(2|g|)\}$, then
\begin{align*} \frac{\Pi(x;g)}{A(g)\pi(x)} &\ge 1-K\left(\frac{\log\log{x}}{\log{x}} + \frac{\log\log{2|g|}}{\log{x}}\right) \\
&\ge 1-\frac{K}{B} - K\frac{\log\log{x}}{\log{x}}.
\end{align*}
The right-hand side is positive for large enough $x$, say $x \ge x_1 = x_1(B,K)$, where $x_1$ is a constant that can be assumed to exceed $x_0$. It follows that $p_g \le \max\{x_1, \log^{B}(2|g|)\}$, giving \eqref{eq:pgupper}.

In our next theorem, we pinpoint a numerically explicit value of $B$.

\begin{thm}[assuming GRH]\label{thm:pgupper} The upper bound $p_g \ll \log^{B}(2|g|)$ holds with $B=19$.
\end{thm}

Usually $p_g$ is quite small. For instance, $p_g = 2$ whenever $g$ is odd, while for even $g$, one has $p_g = 3$ one-third of the time (whenever $3\mid g+1$). Proceeding more generally, there are $\phi(p-1)$ primitive roots modulo the prime $p$. So by the Chinese remainder theorem, for each fixed $p$ a random $g$ satisfies $p_g > p$ with probability $\prod_{r \le p} (1-\frac{\phi(r-1)}{r})$. To make the term ``probability'' here rigorous, we can interpret it as limiting frequency, with $g$ sampled from integers satisfying $|g| \le x$, where $x\to\infty$.

This probabilistic viewpoint suggests a reasonable guess for the maximum size of $p_g$ when $|g|\le x$. While $\phi(r-1)/r$ fluctuates as the prime $r$ varies, for the sake of estimating the above product on $r$, we can treat the terms $1-\frac{\phi(r-1)}{r}$ as constant. More precisely, there is a certain real number $\varrho> 1$ such that \[ \prod_{r \le r_k} \left(1-\frac{\phi(r-1)}{r}\right) = \varrho^{-(1+o(1))k}\quad\text{as}\quad k\to\infty,\] where $r_k$ denotes the $k$th prime in the usual order. (We prove this estimate as Lemma \ref{lem:vaughan} below.)
Hence, one might guess that for a given $k$ and $x$, the number of $g$, $|g|\le x$,  with $p_g > r_k$ is $\approx 2x \varrho^{-k}$. (Here $2x$ approximates the size of the sample space of $g$ values.) The expression $2x \varrho^{-k}$ is smaller than $1$ once $k > k_0(x):=\frac{\log{2x}}{\log{\varrho}}$. It is therefore tempting to conjecture that $\max_{|g| \le x} p_g$ is never more than about $p_{k_0(x)}$. (This argument is purely heuristic; it requires ``pretending'' that our probabilities, which were given rigorous meaning only when fixing $k$ and sending $x$ to infinity, can be interpreted uniformly in $k$ and $x$.) This cannot be quite right, as $p_g=\infty$ for even square values of $g$ !
Nevertheless, it seems sensible to guess that $p_g \ll (\log{2|g|})(\log\log{2|g|})$ for all $g\in \Gg$. If correct, this is sharp: In \cite{PS09}, Pomerance and Shparlinski report a a construction of Soundararajan yielding an infinite sequence of positive integers $g$ that (a) are all products of two distinct primes and (b) are squares modulo every odd prime $p\le 0.7(\log{g})(\log\log{g})$.\footnote{Here $0.7$ can be replaced with any constant smaller than $1/\log{4}$.} These $g$ satisfy $p_{4g} \gg \log{(4g)}\log\log{(4g)}$.

This same perspective suggests that the ``probability'' that a random integer $g$ satisfies $p_g = p$ is given by
\begin{equation}\label{eq:deltapdef} \delta_p:= \frac{\phi(p-1)}{p} \prod_{r < p} \left(1-\frac{\phi(r-1)}{r}\right).\end{equation}
Taking this for granted and proceeding formally, $\E[p_g] = \sum_{p} p \delta_p$. Using Theorem \ref{thm:pgupper}
, we give a GRH-conditional proof that this sum represents the honest average of $p_g$ over $g \in \Gg$. 

\begin{cor}\label{cor:pgaverage} We have that $\sum_{p} p \delta_p < \infty$. Furthermore, assuming GRH,
\begin{equation}\label{eq:pgaverage} \lim_{x\to\infty} \frac{1}{2x} \sum_{g \in \Gg,\, |g|\le x} p_g = \sum_{p} p \delta_p. \end{equation}
Here $\delta_p$ is as defined in \eqref{eq:deltapdef}.
\end{cor}
\noindent (We divide by $2x$, as there are $2x+O(x^{1/2})$ integers $g\in \Gg$ with $|g| \le x$.) There seems no hope at present of proving Corollary \ref{cor:pgaverage} unconditionally: If $p_g=\infty$ for even a single value of $g\in \Gg$, then the average becomes meaningless, and we know of no way to rule this out. Infinite values of $p_g$ are not the only enemy: Having $p_g > x\log{x}$ for some $g\in \Gg$, $|g|\le x$ (along a sequence of $x$ tending to infinity) is enough to doom \eqref{eq:pgaverage}.

In an attempt to salvage the situation, one might tamp down the large values of $p_g$ by averaging $\min\{p_g, \psi(x)\}$ for a threshold function $\psi$. In our final theorem on $p_g$, established in \S\ref{sec:tamedavg}, we show that this strategy succeeds for $\psi(x) = x^{\eta}$, for any positive $\eta < \frac12$.
 
\begin{thm}\label{thm:pgavg2} Fix a positive real number $\eta < \frac{1}{2}$. Then
\[ \lim_{x\to\infty} \frac{1}{2x} \sum_{g \in \Gg,\,|g| \le x} \min\{p_g, x^{\eta}\} = \sum_{p} p \delta_p. \]
\end{thm}

Theorem \ref{thm:pgavg2} implies that any estimate of the shape $\max\{p_g: |g| \le x, g \in \Gg\}  \ll x^{\frac{1}{2}-\epsilon}$ would suffice to establish \eqref{eq:pgaverage}.


It would be interesting to prove Theorem \ref{thm:pgavg2} with a less stringent condition on $\eta$, such as $\eta < 1$. But a substantial new idea  seems required to take $\eta$ past $1/2$. As we demonstrate in \S\ref{sec:pgstar}, the problem becomes easier if we look instead at \textsf{almost-primitive roots}, by which we mean the integers $g$ which generate a subgroup of index at most two inside $(\Z/p\Z)^{\times}$.

The problems we have taken up about $p_g$ are dual to those classically considered for $g_p$, the least primitive root modulo the prime $p$. Burgess \cite{burgess62} and Wang \cite{wang59} have shown unconditionally that $g_p \ll p^{\frac14+\epsilon}$ for all primes $p$, while Shoup \cite{shoup92} (sharpening an earlier, qualitatively similar result of Wang, op.\ cit.) has proved under GRH that $g_p \ll r^4 (1+\log{r})^4 \log^{2}{p}$, where $r= \omega(p-1)$. Shoup's upper bound is of size $\log^{2+o(1)}{p}$ for most primes $p$ and is always $O(\log^6 p)$. These pointwise results are stronger than those known for $p_g$, but the story for average values is different. While $g_p$ is conjectured to have a finite, limiting mean value as $p$ varies (among primes sampled in increasing order), this has not been established even assuming GRH. In fact, GRH has so far not yielded a stronger upper bound for $\pi(x)^{-1} \sum_{p\le x} g_p$ than $(\log{x}) (\log\log{x})^{1+o(1)}$ (as $x\to\infty$). This last estimate is due to Elliott and Murata \cite{EM97}. In \S4 of the same paper, Elliott and Murata propose a precise value for the average of $g_p$. Their theoretical expression is rather unwieldy and not easy to compute with. However, extensive direct computations of $g_p$ by Andrzej Paszkiewicz (reported on in \cite{EM97})  suggest $g_p$ has mean value $\approx 4.924$.

\section{Notation}\label{sec:notation}
We use standard notation for arithmetic functions throughout the paper.
In particular, $\mu$ is the M\"{o}bius function, $\Lambda$ is the von Mangoldt function, $\varphi$ is the Euler totient function, and $\omega$ is the prime omega function, which, when evaluated at a nonzero integer $n$, returns the number of distinct prime factors of $n$.
We write $(\cdot/\cdot)$ for the Kronecker symbol; often the ``denominator'' will be a prime $p$, in which case $(\cdot/p)$ may be viewed as a Legendre symbol.

Throughout, the letters $x,y,z,\delta,\epsilon,\eta,\theta,\rho,K,Q,X,Y$ represent positive real variables, the letters $d,e,f,g,h,k,m,n,M,N$ stand for integer variables, and the letters $p,q,r,\ell$ are reserved for prime variables. The integer part of a real number $x$, which is defined as the greatest integer not exceeding $x$, will be denoted by $\lfloor x\rfloor$. We write $\gcd(m,n)$, or sometimes simply $(m,n)$, for the greatest common divisor of $m$ and $n$. For a positive integer $n$, we denote by $P^+(n)$ the largest prime factor of $n$, with the convention that $P^+(1)=1$, and by $P^-(n)$ the least prime factor of $n$, with the convention that $P^-(1)=\infty$.

We will also adopt the standard Landau--Vinogradov asymptotic notation  such as $O$, $o$, $\ll$ and $\gg$, as well as the notation $\sim$ from calculus. Given real-valued functions $X,Y$ of a variable $t$ in a certain range, the relations $X=O(Y)$ and $X\ll Y$ will be used interchangeably to mean that there exists a constant $C>0$ such that $|X|\le CY$ for all $t$ in the considered range. Next, the relation $X\gg Y$ is equivalent to $Y=O(X)$, and the relation $X=o(Y)$ is interpreted as $X/Y\to0$ as $t\to\infty$. And as usual, we write $X\sim Y$ whenever $X/Y\to1$ as $t\to\infty$.

When it comes to prime counting, we denote by $\pi(x)$ the number of primes $p\le x$, and by $\pi(x;d,a)$ the number of primes $p\le x$ satisfying the congruence $p\equiv a\pmod*{d}$. The prime number theorem then states that $\pi(x)$ is well approximated by the logarithmic integral $\Li(x)\colonequals\int_{2}^{x}1/\log t\,\mathrm{d}t$, which is itself asymptotically equivalent to $x/\log x$ as $x\to\infty$. Finally, for any subset $A\subseteq\Z$ the indicator function $1_{A}$ of $A$ is defined by $1_{A}(n)=1$ if $n\in A$ and $1_{A}(n)=0$ otherwise. Analogously, we define, for any logic statement $P$, $1_{P}=1$ if $P$ is true and $1_{P}=0$ if $P$ is false.

\section{A uniform variant of Hooley's formula: Proof of Theorem \ref{thm:uniformasymptotic}}\label{sec:uniformhooley}
The following lemma encodes the input of GRH to the proof. It will be of vital importance both in this section and the next.

\begin{lem}[assuming GRH]\label{lem:GRHlem} Let $g$ be a nonzero integer. For each real number $x\ge 2$ and each $d\in \N$, the count of primes $p\le x$ for which
\begin{equation}\label{eq:failsdtest} p \equiv 1\pmod{d} \quad\text{and}\quad g^{(p-1)/d}\equiv 1\pmod{p}\end{equation}
is 
\[ \frac{\pi(x)}{[\Q(\zeta_d,\sqrt[d]{g}):\Q]} + O(x^{1/2}\log(|g|dx)). \]
Here the implied constant is absolute.
\end{lem}
\begin{proof} Apart from making explicit the dependence on $g$, this result is well-known and present already in \cite{hooley67}. Since dependence on $g$ is crucial for our purposes, we sketch a proof. We first throw out primes dividing $dg$; there are only $O(\log{(|g|d)})$ of these, a quantity  subsumed by our error term. For the remaining primes $p$,
\begin{align*} \eqref{eq:failsdtest} \text{ holds } &\Longleftrightarrow \text{$x^d-g$ has $d$ distinct roots over $\F_p$}\\ &\Longleftrightarrow \text{$x^d-g$ factors over $\F_p$ into $d$ distinct monic linear polynomials} \\
& \Longleftrightarrow p\text{ splits completely in $\Q(\zeta_d,\sqrt[d]{g})$}. \end{align*}
To count primes up to $x$ satisfying this last condition, we apply the GRH-conditional Chebotarev density theorem in the form ($20_{R}$) of \cite{serre81} (in the notation of \cite{serre81}, take $K=\Q$, $E=\Q(\zeta_d,\sqrt[d]{g})$, $C=\{\text{id}\}$, and keep in mind that all primes ramifying in $E$ divide $gd$).
\end{proof}

We now turn to the proof of Theorem \ref{thm:uniformasymptotic}. We follow  Hooley's strategy, but keep a more watchful eye on $g$-dependence in the error terms.

Let $p$ be a prime not dividing $g$. For each prime number $\ell$, we say that $p$ \textsf{fails the $\ell$-test} if 
\[ p\equiv1\pmod{\ell}\quad\text{and}\quad g^{(p-1)/\ell}\equiv 1\pmod{p}; \]
otherwise, we say $p$ \textsf{passes the $\ell$-test}. Then $g$ is a primitive root modulo $p$ precisely when $p$ passes the $\ell$-test for all primes $\ell$. In particular, if we define
\[ \Pi_0(x;g) = \#\{p\le x: p\nmid g, p\text{ passes all $\ell$-tests for $\ell \le\log{x}$}\},\]
then \[ \Pi(x;g) \le \Pi_0(x;g).\] For each squarefree $d\in \N$, let $N_d$ denote the count of primes $p \le x$ which fail the $\ell$-test for each prime $\ell\mid d$. These are precisely the primes $p\le x$ for which \eqref{eq:failsdtest} holds, so that by Lemma \ref{lem:GRHlem} and inclusion-exclusion,
\begin{align}\notag \Pi_0(x;g) &= \sum_{d:\, P^{+}(d) \le \log{x}} \mu(d) N_d \\
&= \pi(x) \sum_{d:\, P^{+}(d) \le \log{x}} \frac{\mu(d)}{[\Q(\zeta_d,\sqrt[d]{g}):\Q]} +O\bigg(x^{1/2} \sum_{d:\, P^{+}(d) \le \log{x}} \mu(d)^2 \log(|g|dx)\bigg).\label{eq:pi00} \end{align}
The error term is readily handled: Each squarefree $d$ with $P^{+}(d)\le \log{x}$ satisfies $d \le \prod_{r \le \log{x}} r \le x^2$, and there are $2^{\pi(\log{x})} = \exp(O(\log{x}/\log\log{x}))$ such values of $d$. Hence, 
\begin{equation}\label{eq:pi01} x^{1/2} \sum_{d:\, P^{+}(d) \le \log{x}} \mu(d)^2 \log(|g|dx) \ll x^{1/2} \log(|g|x) \cdot \exp(O(\log{x}/\log\log{x})) \ll x^{3/5}\log{|g|}. \end{equation}

Turning to the main term, we extract from \cite[pp. 213--214]{hooley67} or \cite[Proposition 4.1]{wagstaff82} that for each squarefree $d \in \N$,
\begin{align}
 [\Q(\zeta_d,\sqrt[d]{g}):\Q] = \frac{d\phi(d)}{\epsilon(d) \gcd(d,h)}, \quad\text{where}\quad \epsilon(d) = \begin{cases}
 2 &\text{if $2g_1 \mid d$ and $g_1\equiv 1\pmod*{4}$},\\
 1 &\text{otherwise}.
 \end{cases}\label{eq:hooleydegrees}
\end{align}
(Actually, what Hooley computes in \cite{hooley67} is the degree of $\Q(\zeta_{d}, \sqrt[d_1]{g})$, where $d_1:= d/\gcd(d,h)$. But this is the same field as $\Q(\zeta_d,\sqrt[d]{g})$, by Kummer theory, since the classes of $g$ and $g^{\gcd(d,h)}$ generate the same subgroup of $\Q(\zeta_d)^{\times}/
(\Q(\zeta_d)^{\times})^d$.) From this, Hooley deduces in \cite{hooley67} that
\begin{equation}\label{eq:hooleyfullsum} \sum_{d} \frac{\mu(d)}{[\Q(\zeta_d,\sqrt[d]{g}):\Q]} = A(g), \end{equation}
where the sum is over all $d \in \N$. We would like to plug this result into \eqref{eq:pi00}, but the corresponding sum in \eqref{eq:pi00} is restricted to $(\log{x})$-smooth values of $d$.  The next lemma estimates the error incurred by replacing the sum over all $d$ by the sum over $(\log{x})$-smooth $d$. 
\begin{lem}\label{lem:intermediate}We have
\begin{equation} \sum_{d:\, P^{+}(d) > \log{x}} \frac{\mu(d)}{[\Q(\zeta_d,\sqrt[d]{g}):\Q]} \ll \frac{\phi(h)}{h} \cdot \frac{\log\log{2|g|}}{\log{x}}. \label{eq:lemintermediate}\end{equation}
\end{lem}
\begin{proof} If $g_1\not\equiv 1\pmod{4}$, then
\begin{align} \sum_{d:\, P^{+}(d) > \log{x}} \frac{\mu(d)}{[\Q(\zeta_d,\sqrt[d]{g}):\Q]}&= \sum_{d:\, P^{+}(d) > \log{x}} \mu(d) \frac{(d,h)}{d\phi(d)} = -\sum_{\ell > \log{x}}\frac{(\ell,h)}{\ell \phi(\ell)} \sum_{d:\,P^{+}(d) < \ell} \mu(d)\frac{(d,h)}{d\phi(d)}\notag \\
&=- \sum_{\ell > \log{x}} \frac{(\ell,h)}{\ell (\ell-1)} \prod_{\substack{r < \ell \\ r\nmid h}}\left(1-\frac{1}{r(r-1)}\right) \prod_{\substack{r <\ell \\ r\mid h}}\left(1-\frac{1}{r-1}\right)\notag\\
&\ll \sum_{\ell > \log{x}} \frac{(\ell,h)}{\ell(\ell-1)} \frac{\phi(h)}{h} \prod_{\substack{r\mid h \\ r \ge \ell}}\left(1+\frac{1}{r}\right).\label{eq:needthislater0} \end{align}
Each $r$ appearing in this last expression has $r > \log{x}$. Furthermore,
\begin{equation}\label{eq:smallproduct} \prod_{\substack{r\mid h \\ r > \log{x}}}\left(1+\frac{1}{r}\right) \le \exp\bigg(\sum_{\substack{r\mid h \\ r >\log{x}}} \frac{1}{r}\bigg) \le \exp\bigg(\frac{1}{\log{x}} \sum_{\substack{r \mid h \\ r > \log{x}}}1\bigg) \le \exp\bigg(\frac{\log{h}}{\log x\cdot \log\log{x}}\bigg) \ll 1, \end{equation}
noting that
\begin{equation}\label{eq:tonote} h \le \frac{\log{|g|}}{\log{2}} < \log^3({2|g|}) \le x  \end{equation}
in the last step. Hence, $\prod_{r\mid h,~r\ge \ell} (1+1/r) \ll 1$, and
\begin{align} \sum_{\ell > \log{x}} \frac{(\ell,h)}{\ell(\ell-1)} \frac{\phi(h)}{h} \prod_{\substack{r\mid h \\ r \ge \ell}}\left(1+\frac{1}{r}\right) &\ll \frac{\phi(h)}{h} \left(\sum_{\substack{\ell > \log{x} \\ \ell \mid h}}\frac{1}{\ell} +\sum_{\substack{\ell > \log{x} \\ \ell \nmid h}}\frac{1}{\ell^2}\right)\notag \\
&\ll \frac{\phi(h)}{h}\left(\frac{1}{\log{x}} \frac{\log{h}}{\log\log{x}} + \frac{1}{\log{x}}\right)\notag\\
&\ll \frac{\phi(h)}{h} \cdot \frac{\log\log{2|g|}}{\log{x}},\label{eq:needthislater}\end{align}
where we take from \eqref{eq:tonote} that $\log{h} \ll \log\log{2|g|}$. The assertion of Lemma \ref{lem:intermediate} now follows from \eqref{eq:needthislater0} and \eqref{eq:needthislater}, when $g_1 \not\equiv 1\pmod{4}$.

When $g_1\equiv 1\pmod{4}$, the argument is similar, but the details are slightly more involved. In this case,
\[ \sum_{d:\, P^{+}(d) > \log{x}} \frac{\mu(d)}{[\Q(\zeta_d,\sqrt[d]{g}):\Q]} = 
\sum_{d:\, P^{+}(d) > \log{x}} \mu(d) \frac{(d,h)}{d\phi(d)}  + \sum_{\substack{d:\, P^{+}(d) > \log{x} \\ 2g_1 \mid d}} \mu(d) \frac{(d,h)}{d\phi(d)}.\]
The first right-hand sum appeared earlier and was shown to be $O(\frac{\phi(h)}{h} \frac{\log\log{2|g|}}{\log{x}})$ (see \eqref{eq:needthislater0} and following). To finish off the lemma, it suffices to show that the second right-hand sum is bounded by this same $O$-expression. We rewrite
\begin{equation}\label{eq:complicatedcase}
\sum_{\substack{d:\, P^{+}(d) > \log{x} \\ 2g_1 \mid d}} \mu(d) \frac{(d,h)}{d\phi(d)} = -\sum_{\ell >\log{x}} \frac{(\ell,h)}{\ell \phi(\ell)} \sum_{\substack{d:\, P^{+}(d) < \ell \\ 2g_1 \mid \ell d}} \mu(d) \frac{(d,h)}{d\phi(d)}. 
\end{equation}
The right-hand sum on $d$ is empty if $2g_1/(2g_1,\ell)$ has a prime factor $p$ at least $\ell$. Indeed, in that case the condition $2g_1 \mid \ell d$ forces $p\mid d$, contradicting $P^{+}(d) < \ell$.
In all other cases, letting $r$ denote a prime number,
\[ \sum_{\substack{d:\, P^{+}(d) < \ell \\ 2g_1 \mid \ell d}} \mu(d) \frac{(d,h)}{d\phi(d)} = \prod_{r \mid \frac{2g_1}{(2g_1,\ell)}} -\frac{(r,h)}{r(r-1)} \prod_{\substack{r < \ell\\ r \nmid \frac{2g_1}{(2g_1,\ell)}}}\left(1-\frac{(r,h)}{r(r-1)}\right). \]
Keeping in mind that $h$ is odd, we observe that $\frac{(r,h)}{r(r-1)} \le \frac{1}{2}$ for each prime $r$, so that $\frac{(r,h)}{r(r-1)} \le 1- \frac{(r,h)}{r(r-1)}$. Therefore,
\[ \left|\sum_{\substack{d:\, P^{+}(d) < \ell \\ 2g_1 \mid \ell d}} \mu(d) \frac{(d,h)}{d\phi(d)}\right| \le \prod_{r < \ell} \left(1-\frac{(r,h)}{r(r-1)}\right) \le \prod_{\substack{r < \ell \\ r\mid h}} \left(1-\frac{1}{r-1}\right)\\ \ll \frac{\phi(h)}{h} \prod_{\substack{r\mid h \\ r \ge \ell}}\left(1+\frac{1}{r}\right), \]
and referring back to \eqref{eq:complicatedcase},
\[ \sum_{\substack{d:\, P^{+}(d) > \log{x} \\ 2g_1 \mid d}} \mu(d) \frac{(d,h)}{d\phi(d)} \ll \sum_{\ell > \log{x}} \frac{(\ell,h)}{\ell(\ell-1)} \frac{\phi(h)}{h} \prod_{\substack{r\mid h \\ r \ge \ell}}\left(1+\frac{1}{r}\right). \]
To conclude, recall that the right-hand side was estimated as $O(\frac{\phi(h)}{h} \frac{\log\log{2|g|}}{\log{x}})$ already in \eqref{eq:needthislater}.\end{proof}

From \eqref{eq:pi00} and \eqref{eq:pi01}, we have $\Pi_0(x;g) = \pi(x) \left(\sum_{d} \mu(d) [\Q(\zeta_d,\sqrt[d]{g}):\Q]^{-1}\right) + O(x^{3/5}\log|g|)$. Using \eqref{eq:hooleyfullsum} and \eqref{eq:lemintermediate} to handle the sum on $d$, we arrive at the estimate
\begin{equation}\label{eq:multform0} \Pi_0(x;g) = A(g) \pi(x) + O\left(\frac{\phi(h)}{h}\frac{\log\log{2|g|}}{\log x} \pi(x) + x^{3/5}\log|g|\right). \end{equation}

Our next lemma puts the error term in ``multiplicative form''.
\begin{lem}\label{lem:multform} We have
\begin{equation}\label{eq:pi0stepping} \Pi_0(x;g) = A(g)\pi(x) \left(1 + O\left(\frac{\log\log{2|g|}}{\log{x}}\right)\right). \end{equation}
\end{lem}

\begin{proof} Notice that $A_0(g)$, as defined in \eqref{eq:A0eq}, satisfies $A_0(g) \asymp \phi(h)/h$. Recalling the definition \eqref{eq:A1eq}  of $A_1(g)$ in the case when $g_1\equiv 1\pmod{4}$, we see that the subtracted term in \eqref{eq:A1eq} always has absolute value at most $1$. In fact, that absolute value is at most $1/3$ unless $g_1=-3$, in which case $\mu(|g_1|) = -1$. Hence, $\frac{2}{3} \le A_1(g) \le 2$, and 
\[ A(g) = A_0(g) A_1(g) \asymp \frac{\phi(h)}{h}. \]
Furthermore, $h < x$ (see \eqref{eq:tonote}), so that
\[ \frac{h}{\phi(h)} \ll \log\log{3h} \ll \log\log{x}, \]
while (again from \eqref{eq:tonote})
\[ \log|g| \ll x^{1/3} = x^{3/8}/x^{1/24}. \]
Therefore,
\begin{align*} \frac{\phi(h)}{h}\pi(x)\cdot  \frac{\log\log{2|g|}}{\log x} + x^{3/5}\log|g|
&\ll A(g) \pi(x) \left(\frac{\log\log{2|g|}}{\log{x}}+ \frac{(h/\phi(h))\log|g|}{x^{3/8}}\right) \\
&\ll A(g) \pi(x) \left(\frac{\log\log{2|g|}}{\log{x}}+ \frac{\log\log{x}}{x^{1/24}} \right) \\
&\ll A(g) \pi(x) \frac{\log\log{2|g|}}{\log{x}}.\end{align*}
The assertion \eqref{eq:pi0stepping} of Lemma \ref{lem:multform} now follows from \eqref{eq:multform0}.
\end{proof}

Next, we investigate the difference $\Pi_0(x;g) - \Pi(x;g)$. If the prime $p \le x$ is counted by $\Pi_0(x;g)$ but not $\Pi(x;g)$, then $p$ passes the $\ell$-tests for all $\ell \le \log{x}$ but fails the $\ell$-test for some $\ell > \log{x}$. Set
\[ x_1 = \log{x}, \quad x_2 = x^{1/2}(\log{x})^{-2} (\log{|g|})^{-1}, \quad x_3 = x^{1/2} (\log{x})^2 \log{|g|},\]
and put
\[ I_1 = (x_1, x_2], \quad I_2 = (x_2, x_3], \quad I_3 = (x_3,\infty). \]
For $j\in\{1,2,3\}$, let $E_j$ denote the count of primes $p\le x$, $p\nmid g$, which fail the $\ell$-test for the first time for an $\ell \in I_j$. Then
\begin{equation}\label{eq:pixglower} \Pi_0(x;g) \ge \Pi(x;g) \ge \Pi_0(x;g) - E_1 - E_2 - E_3.\end{equation}
We proceed to estimate the $E_j$ in turn.

\begin{lem} We have
\begin{equation}\label{eq:A1suitable} E_1 \ll A(g) \pi(x) \left(\frac{\log\log{2|g|}}{\log{x}} + \frac{\log\log{x}}{\log{x}}\right).\end{equation}
\end{lem}

\begin{proof} Recall that for a prime $\ell$, we are using $N_{\ell}$ for the number of primes $p\le x$, $p\equiv 1\pmod{\ell}$, for which $g^{(p-1)/\ell}\equiv 1\pmod{p}$. Invoking Lemma \ref{lem:GRHlem}, and keeping in mind that $[\Q(\zeta_\ell,\sqrt[\ell]{g}):\Q] \gg \ell^2/(\ell,h)$ by \eqref{eq:hooleydegrees}, we find that
\begin{align*} E_1 \le \sum_{\ell \in I_1} N_{\ell} &\ll \sum_{\ell \in I_1} \left(\pi(x) \frac{(\ell,h)}{\ell^2} + x^{1/2}\log(|g|\ell x)\right) \\
&\ll \pi(x) \Bigg(\sum_{\ell >\log{x}} \frac{1}{\ell^2} + \sum_{\substack{\ell > \log{x} \\ \ell \mid h}} \frac{1}{\ell}\Bigg) + x^{1/2}\log(|g|x) \cdot \pi(x_2).
\end{align*}
Since $h < x$ and $\log h \ll \log\log{2|g|}$,
\begin{multline*} \sum_{\ell >\log{x}} \frac{1}{\ell^2} + \sum_{\substack{\ell > \log{x} \\ \ell \mid h}} \frac{1}{\ell} \ll \frac{1}{\log{x}\cdot\log\log{x}} + \frac{1}{\log{x}} \frac{\log{h}}{\log\log{x}} \ll \frac{\log\log{2|g|}}{\log{x}\cdot\log\log{x}} \\ = \frac{\phi(h)}{h} \left(\frac{h/\phi(h)}{\log x} \frac{\log\log{2|g|}}{\log\log{x}}\right) \ll \frac{\phi(h)}{h} \left(\frac{\log\log{x}}{\log x} \frac{\log\log{2|g|}}{\log\log{x}}\right)  = \frac{\phi(h)}{h} \frac{\log\log{2|g|}}{\log{x}}, \end{multline*}
so that
\[ \pi(x) \Bigg(\sum_{\ell >\log{x}} \frac{1}{\ell^2} + \sum_{\substack{\ell > \log{x} \\ \ell \mid h}} \frac{1}{\ell}\Bigg) \ll A(g)\pi(x) \cdot \frac{\log\log{2|g|}}{\log{x}}. \]
We are assuming that $x \ge (\log{2|g|})^3$. Hence, \begin{equation}\label{eq:lowerx2bound} x_2\ge x^{1/6} (\log{x})^{-2} > x^{1/7} \end{equation} for all $x$ exceeding a certain absolute constant, and $\log{x_2} \gg \log{x}$. Thus, $\pi(x_2) \ll x_2 (\log{x})^{-1} = x^{1/2} (\log{x})^{-3} (\log{|g|})^{-1}$, and 
\begin{align*} x^{1/2}\log(|g| x) \cdot \pi(x_2) &\ll \frac{x}{(\log{x})^3 \log{|g|}} (\log{|g|x}) \ll \pi(x) \frac{\log{|g|x}}{(\log{x})^2\log|g|} \\ &\ll \frac{\pi(x)}{\log{x}} = \frac{\phi(h)}{h} \pi(x) \cdot \frac{h/\phi(h)}{\log{x}}  \ll A(g) \pi(x)  \frac{\log\log{x}}{\log{x}}. \end{align*}
Collecting our observations yields \eqref{eq:A1suitable}.\end{proof}

\begin{lem} We have
\begin{equation}\label{eq:A2suitable} E_2 \ll \pi(x) A(g) \left(\frac{\log\log{x}}{\log{x}}+ \frac{\log\log{2|g|}}{\log{x}} \right).\end{equation}
\end{lem}

\begin{proof} Let $\ell$ be a prime dividing $h$. Then every prime $p\equiv 1\pmod{\ell}$, with $p$ not dividing $g$, satisfies 
\[ g^{(p-1)/\ell} \equiv 1\pmod{p}, \]
as $g$ is an $\ell$th power. Hence, in order for a prime $p$ (not dividing $g$) to pass the $\ell$-test, it must be that $p\not\equiv 1\pmod{\ell}$. By assumption, the primes counted in $E_2$ pass the $\ell$-test for all $\ell \le x_2$, and hence for all $\ell \le x^{1/7}$ (see \eqref{eq:lowerx2bound}). So if we let $h'$ denote the $x^{1/7}$-smooth part of $h$, then each prime $p$ counted in $E_2$ has $(p-1,h')=1$. Since $p$ also fails the $\ell$-test for some $\ell \in I_2$, 
\[ E_2 \le \sum_{\ell \in I_2} \sum_{\substack{p \le x \\ (p-1,h')=1 \\ p\equiv 1\pmod*{\ell}}} 1. \]
Each prime $p$ counted by the inner sum has the form $p = 1 + \ell m$. Here $0 < m < x/\ell$, and $m$ avoids the residue classes $0\bmod{r}$ for all primes $r\mid h$, $r \le x^{1/7}$, as well as the residue classes of $-1/\ell\bmod{r}$ for each prime $r < \ell$. Moreover, for each $\ell \in I_2$, we have $\ell > x_2 > x^{1/7}$ as well as $x/\ell \ge x/x_3 = x_2 > x^{1/7}$. Applying Brun's sieve, 
\[ \sum_{\substack{p \le x \\ (p-1,h')=1 \\ p\equiv 1\pmod*{\ell}}} 1 \ll \frac{x}{\ell} \prod_{r \le x^{1/7}} \left(1-\frac{1 + 1_{r\mid h}}{r}\right) \ll \frac{x}{\ell \log x} \prod_{\substack{r \le x^{1/7} \\ r \mid h}}\left(1-\frac{1}{r}\right) \ll \frac{\pi(x)}{\ell} \frac{\phi(h)}{h} \prod_{\substack{r>x^{1/7} \\ r \mid h}}\left(1+\frac{1}{r}\right).\]
(Here and below, ``Brun's sieve'' can be taken to refer to Theorem 2.2 on p.\ 68 of \cite{HR74}.) We have from \eqref{eq:smallproduct} that the final product on $r$ is $O(1)$. Thus,
\[ \sum_{\ell \in I_2} \sum_{\substack{p \le x \\ (p-1,h')=1 \\ p\equiv 1\pmod*{\ell}}} 1 \ll \pi(x) \frac{\phi(h)}{h} \sum_{\ell \in I_2} \frac{1}{\ell} \\
\ll \pi(x) \frac{\phi(h)}{h} \left(\frac{\log\log{x}}{\log{x}}+ \frac{\log\log{2|g|}}{\log{x}} \right), \]
using Mertens' theorem \cite[Theorem 2.7(d)]{MV06} to estimate the sum on $\ell$. Recalling that $A(g) \asymp \phi(h)/h$, we obtain
\eqref{eq:A2suitable}. \end{proof}

\begin{lem} We have
\[ E_3 \ll A(g)\pi(x) \cdot \frac{\log\log{x}}{\log{x}}.\]
\end{lem} 

\begin{proof} Each $p$ counted in $E_3$ has $g^{(p-1)/\ell}\equiv 1\pmod{p}$ for some $\ell > x_3$. Thus, the multiplicative order of $g$ mod $p$ is smaller than $x/x_3=x_2$, and $p$ divides $g^{m}-1$ for some natural number $m < x_2$. The number of distinct prime factors of $g^m-1$ is $O(m\log|g|)$, and so 
\[ E_3 \ll \log|g| \sum_{m < x_2} m \ll x_2^2 \log|g| = \frac{x}{(\log^4{x}) (\log|g|)}. \]
In particular,
\begin{equation}\label{eq:A3suitable} E_3 \ll \frac{\pi(x)}{\log{x}} = \frac{\phi(h)}{h}\pi(x) \cdot \frac{h/\phi(h)}{\log{x}}  \ll A(g)\pi(x) \cdot \frac{\log\log{x}}{\log{x}}, \end{equation}as desired.\end{proof}

Combining \eqref{eq:pi0stepping}, \eqref{eq:pixglower},  \eqref{eq:A1suitable}, \eqref{eq:A2suitable}, and \eqref{eq:A3suitable}, 
\[ \Pi(x;g) = A(g) \pi(x) \left(1 + O\left(\frac{\log\log{x}}{\log{x}} + \frac{\log\log{2|g|}}{\log{x}}\right)\right). \]
This completes the proof of Theorem \ref{thm:uniformasymptotic}.

\section{An explicit upper bound for the least Artin prime $p_g$: Proof of Theorem \ref{thm:pgupper}}\label{sec:explicitp_g}
Now we turn to the proof of Theorem \ref{thm:pgupper}. We may assume that $|g|$ is sufficiently large. Let $x=\log^B|g|$ with $B=19$, and put $W=\prod_{2<p\le\log x}p$. Denote by $\Ss$ the set of primes $p\le x$ with $(g/p)=-1$ and $\gcd(p-1,W)=1$. 

First of all, let us estimate the number of elements in $\Ss$. We observe that
\[\#\Ss=\frac{1}{2}\sum_{\substack{p \le x,\,p\nmid g\\(p-1,W)=1}} \left(1-(g/p)\right)=\frac{1}{2}\sum_{\substack{p \le x\\(p-1,W)=1}} \left(1-(g/p)\right)+O(\omega(g)).\]
By inclusion-exclusion, we have
\begin{align*}
\sum_{\substack{p \le x\\(p-1,W)=1}} \left(1-(g/p)\right)&=\sum_{p \le x} \left(1-(g/p)\right) \sum_{\substack{d\mid p-1 \\ d \mid W}}\mu(d)\\
&=\sum_{d\mid W}\mu(d)\sum_{\substack{p\le x\\p\equiv1\pmod*{d}}}\left(1-(g/p)\right)\\
&=\sum_{d\mid W}\mu(d)\pi(x;d,1)-\sum_{d\mid W}\mu(d)\sum_{\substack{p\le x\\p\equiv1\pmod*{d}}}(g/p),
\end{align*}
where $\pi(x;d,1)$ denotes the number of primes $p\le x$ with $p\equiv1\pmod{d}$. Hence,
\begin{align*}
\#\Ss=\frac{1}{2}\sum_{d\mid W}\mu(d)\pi(x;d,1)-\frac{1}{2}\sum_{d\mid W}\mu(d)\sum_{\substack{p\le x\\p\equiv1\pmod*{d}}}(g/p)+O(\log|g|),
\end{align*}
since $\omega(g)\le2\log|g|$. To estimate the first sum above, we appeal to \cite[Corollary 13.8]{MV06}, the GRH-conditional prime number theorem for primes in arithmetic progressions, to obtain
\begin{align*}
\sum_{d\mid W}\mu(d)\pi(x;d,1)=\sum_{d\mid W}\mu(d)\left(\frac{\Li(x)}{\varphi(d)}+O\left(x^{1/2}\log x\right)\right)&=\tilde{A}_0(g)\Li(x)+O\left(2^{\pi(\log x)}x^{1/2}\log x\right)\\
&=\tilde{A}_0(g)\Li(x)+O\left(x^{1/2+o(1)}\right),
\end{align*}
where 
\[\tilde{A}_0(g)=\sum_{d\mid W}\frac{\mu(d)}{\varphi(d)}=\prod_{2<q\le\log x}\left(1-\frac{1}{q-1}\right).\]
In addition, we can rewrite 
\[\sum_{\substack{p\le x\\p\equiv1\pmod*{d}}}(g/p)=\frac{1}{\varphi(d)}\sum_{\chi\pmod*{d}}\sum_{\substack{p \le x\\ p\nmid g}}\chi(p)(g/p),\]
by the orthogonality relations of Dirichlet characters, where the outer sum on the right-hand side runs over all Dirichlet characters $\chi\pmod{d}$. It follows that
\begin{equation}\label{eq:cardS1}
\#\Ss=\frac{\tilde{A}_0(g)}{2}\Li(x)-\frac{1}{2}\sum_{d\mid W}\frac{\mu(d)}{\varphi(d)}\sum_{\chi\pmod*{d}} \sum_{\substack{p \le x\\ p\nmid g}}\chi(p)(g/p)+O\left(x^{1/2+o(1)}\right).
\end{equation}
To estimate the triple sum in \eqref{eq:cardS1}, we recall that $\Q(\sqrt{g})=\Q(\sqrt{g_1})$, where $g_1\ne 1$ is the unique squarefree integer with $g_1 (\Q^{\times})^2 = g(\Q^{\times})^2$. Let $\Delta$ be the discriminant of $\Q(\sqrt{g_1})$. Then $(g/p) = (\Delta/p)$ for all odd primes $p$ not dividing $g$. For these primes $p$, $\chi(p)(g/p)$ can be viewed as the value at $p$ of a character $\psi_{\chi,g}\pmod{|\Delta| d}$. The character $\psi_{\chi,g}$ is non-principal unless $\chi$ is induced by the primitive character $(\Delta/\cdot)\pmod{|\Delta|}$. For that to occur, one needs $\Delta\mid d$; in that eventuality, to each $d$ there corresponds exactly one character $\chi\pmod{d}$ for which $\psi_{\chi,g}$ is trivial. All of the $d$ appearing above are odd, squarefree, and divide $W$, so in order for $\Delta$ to divide $d$ we need $\Delta$ to be a squarefree divisor of $W$. This forces $\Delta = g_1 \equiv 1\pmod{4}$ and requires that $g_1 \mid W$. By 
\cite[Theorem 13.7]{MV06}, the GRH-conditional estimates for character sums over primes, we have
\begin{align*}
\frac{1}{2}\sum_{d\mid W}\frac{\mu(d)}{\varphi(d)}\sum_{\chi\pmod*{d}} \sum_{\substack{p \le x\\ p\nmid g}}\chi(p)(g/p)&=\frac{1}{2}\sum_{d\mid W}\frac{\mu(d)}{\varphi(d)}\left(1_{g_1\mid d}\cdot1_{4\mid(g_1-1)}\Li(x)+O\left(\varphi(d)x^{1/2}\log(dx)\right)\right)\\
&=\frac{1_{4\mid(g_1-1),\,g_1\mid W}}{2}\Li(x)\sum_{g_1\mid d,\,d\mid W}\frac{\mu(d)}{\varphi(d)}+O\left(2^{\pi(\log x)}x^{1/2}\log x\right)\\
&=\frac{1_{4\mid(g_1-1),\,g_1\mid W}}{2}\cdot\frac{\mu(g_1)}{\varphi(g_1)}\Li(x)\sum_{d\mid (W/g_1)}\frac{\mu(d)}{\varphi(d)}+O\left(x^{1/2+o(1)}\right)\\
&=\frac{1_{4\mid(g_1-1),\,g_1\mid W}}{2}\cdot\frac{\mu(g_1)}{\varphi(g_1)}\Li(x)\prod_{q\mid(W/g_1)}\left(1-\frac{1}{q-1}\right)+O\left(x^{1/2+o(1)}\right)\\
&=\frac{\tilde{A}_0(g)(1-\tilde{A}_1(g))}{2}\Li(x)+O\left(x^{1/2+o(1)}\right),
\end{align*}
where 
\[\tilde{A}_1(g)\colonequals 1-1_{4\mid(g_1-1),\,g_1\mid W}\frac{\mu(g_1)}{\varphi(g_1)}\prod_{q\mid g_1}\left(1-\frac{1}{q-1}\right)^{-1}=1-1_{4\mid(g_1-1),\,g_1\mid W}\prod_{q\mid g_1}\frac{-1}{q-2}.\]
Inserting this estimate in \eqref{eq:cardS1} yields
\begin{equation}\label{eq:cardS2}
\#\Ss=\frac{\tilde{A}_0(g)\tilde{A}_1(g)}{2}\Li(x)+O\left(x^{1/2+o(1)}\right).
\end{equation}
It is worth noting that 
\begin{align*}
\tilde{A}_0(g)=\prod_{2<q\le\log x}\left(1-\frac{1}{q-1}\right)&=\prod_{2<q\le\log x}\left(1-\frac{1}{q}\right)\prod_{2<q\le\log x}\left(1-\frac{1}{q-1}\right)\left(1-\frac{1}{q}\right)^{-1}\\
&=\left(1+O\left(\frac{1}{\log\log x}\right)\right)\frac{2C_2 \mathrm{e}^{-\gamma}}{\log\log x}
\end{align*}
by Mertens' theorem \cite[Theorem 2.7(e)]{MV06}, where $\gamma=0.577215...$ is the Euler–Mascheroni constant, and that
\[\frac{2}{3}=\tilde{A}_1(-15)\le\tilde{A}_1(g)\le \tilde{A}_1(-3)=2,\]
where 
\[C_2\colonequals\prod_{q>2}\left(1-\frac{1}{q-1}\right)\left(1-\frac{1}{q}\right)^{-1}=\prod_{q>2}\left(1-\frac{1}{(q-1)^2}\right)\]
is the twin prime constant. Thus, the main term in \eqref{eq:cardS2} is of order $\Li(x)/\log\log x$.

Next, we estimate the number of $p\in\Ss$ modulo which $g$ is not a primitive root. To this end, we count those $p\in\Ss$ which fail the $\ell$-test for some $\ell>\log x$. Such an $\ell$ falls necessarily into one of the following four intervals:
\begin{alignat*}{2}
&J_1\colonequals(\log x,y_1], \quad\quad\quad &&J_2\colonequals (y_1,y_2],\\
&J_3\colonequals (y_2,x^{\alpha}],\quad\quad\quad &&J_4\colonequals (x^{\alpha},\infty),
\end{alignat*}
where $\alpha\in(10/19,1)$ is fixed, and
\[y_1\colonequals \frac{x^{1/2}}{(\log|g|)\log^2 x},\quad\quad y_2\colonequals x^{1/2-1/\log\log x}.\]
We start with $J_1$. Suppose first that $\ell\nmid h$. Applying Lemma \ref{lem:GRHlem} as in the proof of Theorem \ref{thm:uniformasymptotic}, with the asymptotic relation $\pi(x)\sim \Li(x)$ in mind, we see that the count of $p\in\Ss$ that fail the $\ell$-test for some $\ell \in J_1$ with $\ell\nmid h$ is 
\[ \ll \sum_{\ell \in J_1} \left(\frac{\Li(x)}{\ell^2} + x^{1/2}\log(|g|\ell x)\right)\ll \Li(x)\sum_{\ell>\log x} \frac{1}{\ell^2} + x^{1/2}\pi(y_1)\log(|g|)\ll \frac{\Li(x)}{\log{x}}, \]
which is negligible compared to the main term in \eqref{eq:cardS2}. In the case where $\ell\mid h$, we observe that a prime 
$p \in \Ss$ failing the $\ell$-test satisfies $p\equiv 1\pmod{\ell}$ and $\gcd(p-1,W)=1$. For each $\ell \in J_1$, the number of such $p\le x$ is
\begin{align*}
\sum_{\substack{p\le x\\p\equiv 1\pmod*{\ell}\\(p-1,W)=1}}1\le x^{1/3}+\sum_{\substack{m\le x/\ell\\ (m,W)=1\\P^-(\ell m+1)>x^{1/3}}}1&\ll x^{1/3}+\frac{x}{\ell}\prod_{q\le x^{1/3}}\left(1-\frac{1_{q\mid W}+1_{q\ne \ell}}{q}\right)\\
&\ll x^{1/3}+\frac{x}{\ell}\prod_{q\mid W}\left(1-\frac{1}{q}\right)\prod_{\substack{q\le x^{1/3}\\q\ne \ell}}\left(1-\frac{1}{q}\right)\\
&\ll\frac{\Li(x)}{\ell\log\log x},  
\end{align*}
by Brun's sieve. Summing this on $\ell > \log{x}$ with $\ell \mid h$ gives
\[ \ll \frac{\Li(x)}{\log\log{x}} \sum_{\substack{\ell > \log{x}\\\ell \mid h}} \frac{1}{\ell} \ll \frac{\Li(x)}{(\log x)\log\log{x}} \sum_{\substack{\ell > \log{x}\\\ell \mid h}}1\ll\frac{\Li(x)}{(\log x)\log\log{x}}\cdot\frac{\log h}{\log\log{x}}. \]
Since $h \ll \log{|g|} = x^{1/B}$, this is $\ll \Li(x)/(\log\log x)^2$, which is also negligible compared to the main term in \eqref{eq:cardS2}.

Moving on to $J_2$, we seek to bound the number of primes $p\in\Ss$ failing the $\ell$-test for some $\ell\in J_2$. Such a prime $p$ certainly satisfies $p\le x$, $\gcd(p-1,W)=1$, and $p\equiv1\pmod{\ell}$. Using inclusion-exclusion and invoking \cite[Corollary 13.8]{MV06} again, we find that for each $\ell \in J_2$, the number of such $p$ is 
\begin{align*}
\sum_{\substack{p\le x\\p\equiv 1\pmod*{\ell}\\(p-1,W)=1}}1=\sum_{d\mid W}\mu(d)\pi(x;\ell d,1)&=\sum_{d\mid W}\mu(d)\left(\frac{\Li(x)}{\varphi(\ell d)}+O\left(x^{1/2}\log x\right)\right)\\
&=\frac{\Li(x)}{\varphi(\ell)}\sum_{d\mid W}\frac{\mu(d)}{\varphi(d)}+O\left(2^{\pi(\log x)}x^{1/2}\log x\right)\\
&=\frac{\tilde{A}_0(g)}{\ell-1}\Li(x)+O\left(2^{\pi(\log x)}x^{1/2}\log x\right).
\end{align*}

Summing this on $\ell \in J_2$ shows that the number of primes $p\in\Ss$ failing the $\ell$-test for some $\ell\in J_2$ is
\begin{align*}
&\quad\quad\sum_{\ell\in J_2}\left(\frac{\tilde{A}_0(g)}{\ell-1}\Li(x)+O\left(2^{\pi(\log x)}\sqrt{x}\log x\right)\right)\\
&=\left(\log\frac{\log y_2}{\log y_1}+O\left(\frac{1}{\log y_1}\right)\right)\tilde{A}_0(g)\Li(x)+O\left(2^{\pi(\log x)}\pi(y_2)\sqrt{x}\log x\right)\\
&=\left(\log\frac{B}{B-2}+O\left(\frac{1}{\log\log x}\right)\right)\tilde{A}_0(g)\Li(x)+O\left(x^{1-(1-\log 2+o(1))/\log\log x}\right)\\
&=\left(\log\frac{B}{B-2}+O\left(\frac{1}{\log\log x}\right)\right)\tilde{A}_0(g)\Li(x),
\end{align*}
where we have made use of Mertens' theorem in the first equality and the prime number theorem and the relation $x=\log^B|g|$ in the second equality.

Now we turn to $J_3$. As in the treatment of $J_2$, we shall only use that a prime $p\in\Ss$ failing the $\ell$-test satisfies $p\equiv 1\pmod{\ell}$ and that $\gcd(p-1,W)=1$. However, \cite[Corollary 13.8]{MV06} loses its strength in this case, for most $\ell\in J_3$ go way beyond $x^{1/2}$. To get around this issue, we resort to the following ``arithmetic large sieve'' inequality due to Montgomery (see \cite[Chapter 3]{montgomery71} and \cite[\S 9.4]{FI10}) to obtain an asymptotically explicit upper bound for the number of primes $p\le x$ satisfying $p\equiv 1\pmod{\ell}$ and $\gcd(p-1,W)=1$, rather than pursue an asymptotic formula for this count. 

\begin{ALS} Let $Q\ge1$. To each prime $p \le Q$, associate $\nu(p) < p$ residue classes modulo $p$. For every pair of integers $M,N$, with $N>0$, the number of integers in $[M+1,M+N]$ avoiding the distinguished residue classes mod $p$ for all primes $p \le Q$ is bounded above by 
\[ \frac{N+Q^2}{J}, \quad\text{where}\quad J:= \sum_{n \le Q} \mu^2(n) \prod_{p\mid n} \frac{\nu(p)}{p-\nu(p)}. \] 
\end{ALS}

By the arithmetic large sieve, the count of $p\le x$ corresponding to a given $\ell\in J_3$ is at most
\begin{equation}\label{eq:ALSJ_3}
\sum_{\substack{m\le x/\ell\\ (m,V)=1\\P^-(\ell m+1)>(x/\ell)^{\beta}}}1\le \left(\frac{x}{\ell}+\left(\frac{x}{\ell}\right)^{2\beta}\right)\left(\sum_{n \le (x/\ell)^{\beta}} \mu(n)^2 \prod_{q\mid n}\frac{\nu(q)}{q-\nu(q)}\right)^{-1}   
\end{equation}
where $\beta=\beta(x)=1/2-1/\log\log x$, $V$ is the product of all odd primes not exceeding $\log x/\log\log x$, and $\nu(q)=1_{q\mid V}+1$. Here we have exploited the facts that $V\mid W$ and that $(x/\ell)^{\beta}<\ell$ for every $\ell\in J_3$. To handle the sum on the right-hand side, we observe that $V=x^{(1+o(1))/\log\log x}=(x/\ell)^{O(1/\log\log x)}$ and that
\begin{equation}\label{eq:ALSJ_3a}
\sum_{n \le (x/\ell)^{\beta}} \mu(n)^2 \prod_{q\mid n}\frac{\nu(q)}{q-\nu(q)}\ge\left(\sum_{d\mid V}\mu(d)^2\prod_{q\mid d}\frac{2}{q-2}\right)\left(\sum_{\substack{m \le (x/\ell)^{\beta}/V\\(m,V)=1}}\mu(m)^2\prod_{q\mid m}\frac{1}{q-1}\right).   
\end{equation}
It is easy to see that 
\begin{equation}\label{eq:ALSJ_3b}
\sum_{d\mid V}\mu(d)^2\prod_{q\mid d}\frac{2}{q-2}=\prod_{q\mid V}\left(1+\frac{2}{q-2}\right)=\left(1+O\left(\frac{\log\log\log x}{\log\log x}\right)\right)\prod_{q\mid W}\left(1+\frac{2}{q-2}\right). 
\end{equation}
In addition, we have
\[\sum_{\substack{m \le (x/\ell)^{\beta}/V\\(m,V)=1}}\mu(m)^2\prod_{q\mid m}\frac{1}{q-1}=\sum_{\substack{m \le (x/\ell)^{\beta}/V\\(m,V)=1}}\frac{\mu(m)^2}{\varphi(m)}\ge\frac{\varphi(V)}{V}\sum_{m \le (x/\ell)^{\beta}/V}\frac{\mu(m)^2}{\varphi(m)},\]
where the last inequality follows from
\[\sum_{n\le z}\frac{\mu(n)^2}{\varphi(n)}\le\Bigg(\sum_{d\mid a}\frac{\mu(d)^2}{\varphi(d)}\Bigg)\Bigg(\sum_{\substack{m \le z\\(m,a)=1}}\frac{\mu(m)^2}{\varphi(m)}\Bigg)\]
and 
\[\sum_{d\mid a}\frac{\mu(d)^2}{\varphi(d)}=\frac{a}{\varphi(a)}\]
for all $z\ge1$ and $a\in\N$. Since an application of \cite[eq. (3.18)]{MV06} yields
\[\sum_{m \le (x/\ell)^{\beta}/V}\frac{\mu(m)^2}{\varphi(m)}>\log\frac{(x/\ell)^{\beta}}{V}=\left(\frac{1}{2}+O\left(\frac{1}{\log\log x}\right)\right)\log(x/\ell),\]
we obtain
\begin{align*}
\sum_{\substack{m \le (x/\ell)^{\beta}/V\\(m,V)=1}}\mu(m)^2\prod_{q\mid m}\frac{1}{q-1}&\ge\left(\frac{1}{2}+O\left(\frac{1}{\log\log x}\right)\right)\frac{\varphi(V)}{V}\log(x/\ell)\\
&=\left(\frac{1}{2}+O\left(\frac{\log\log\log x}{\log\log x}\right)\right)\frac{\varphi(W)}{W}\log(x/\ell).
\end{align*}
Inserting this estimate and \eqref{eq:ALSJ_3b} in \eqref{eq:ALSJ_3a} yields
\begin{align*}
\sum_{n \le (x/\ell)^{\beta}} \mu(n)^2 \prod_{q\mid n}\frac{\nu(q)}{q-\nu(q)}&\ge\left(\frac{1}{2}+O\left(\frac{\log\log\log x}{\log\log x}\right)\right)\frac{\varphi(W)}{W}\log(x/\ell)\prod_{q\mid W}\left(1+\frac{2}{q-2}\right)\\
&=\left(\frac{1}{2}+O\left(\frac{\log\log\log x}{\log\log x}\right)\right)\tilde{A}_0(g)^{-1}\log(x/\ell).
\end{align*}
Combining the above with \eqref{eq:ALSJ_3}, we find that the count of $p\le x$ corresponding to a given $\ell\in J_3$ is at most
\[\left(2+O\left(\frac{\log\log\log x}{\log\log x}\right)\right)\tilde{A}_0(g)\frac{x}{\ell\log(x/\ell)}=\left(2+O\left(\frac{\log\log\log x}{\log\log x}\right)\right)\tilde{A}_0(g)\frac{\Li(x)\log x}{\ell\log(x/\ell)}.\]
Summing this on $\ell \in J_3$, we see that the count of $p\le x$ in consideration is at most
\[\left(2+O\left(\frac{\log\log\log x}{\log\log x}\right)\right)\tilde{A}_0(g)\Li(x)\log x\sum_{\ell\in J_3}\frac{1}{\ell\log(x/\ell)}. \]
By Mertens' theorem and partial summation \cite[eq. (A.4), p. 488]{MV06}, we have
\begin{align*}
\sum_{\ell\in J_3}\frac{1}{\ell\log(x/\ell)}&=\int_{J_3}\frac{1}{\log(x/t)}\,\mathrm{d}\left(\sum_{\ell\le t}\frac{1}{\ell}\right)\\
&=\int_{J_3}\frac{\mathrm{d}t}{t(\log t)\log(x/t)}+\int_{J_3}\frac{1}{\log(x/t)}\,\mathrm{d}\left(O\left(\frac{1}{\log t}\right)\right)\\
&=\frac{1}{\log x}\int_{1/2-1/\log\log x}^{\alpha}\frac{\mathrm{d}u}{u(1-u)}+O\left(\frac{1}{(\log x)^2}\right)\\
&=\frac{1}{\log x}\int_{1/2}^{\alpha}\frac{\mathrm{d}u}{u(1-u)}+O\left(\frac{1}{(\log x)\log\log x}\right)\\
&=\left(\log\frac{\alpha}{1-\alpha}+O\left(\frac{1}{\log\log x}\right)\right)\frac{1}{\log x}.
\end{align*}
Hence, the count of $p\le x$ in consideration is at most
\[\left(2\log \frac{\alpha}{1-\alpha}+O\left(\frac{\log\log\log x}{\log\log x}\right)\right)\tilde{A}_0(g)\Li(x).\]

Finally, it remains to estimate the number of primes $p\in\Ss$ failing the $\ell$-test for some $\ell\in J_4$. For each such $p$, the order of $g\bmod{p}$ is smaller than $x^{1-\alpha}$. Thus, $p\mid(g^m-1)$ for some positive integer $m \le x^{1-\alpha}$. The number of distinct prime factors of $g^m-1$ is $O(m\log|g|)$. Hence, the number of primes $p\in\Ss$ failing the $\ell$-test for some $\ell\in J_4$ is at most
\[\sum_{m \le x^{1-\alpha}} m\log|g|\ll x^{2-2\alpha}\log|g|=x^{2-2\alpha+1/B}.\] 
Since $\alpha\in(10/19,1)$, we have $2-2\alpha + 1/B<1$. Thus, $x^{2-2\alpha}\log|g|$ is of smaller order than the main term in \eqref{eq:cardS2}. 

Putting everything together, we deduce that the number of $p\in \Ss$ having $g$ as a primitive root is at least
\[\left(\frac{\tilde{A}_1(g)}{2} - \log\frac{B}{B-2} -2 \log\frac{\alpha}{1-\alpha}+o(1)\right)\tilde{A}_0(g)\Li(x).\]
Since $\tilde{A}_1(g)\ge2/3$, our choice of $B$ guarantees that
\[\frac{\tilde{A}_1(g)}{2} - \log\frac{B}{B-2} -2 \log\frac{\alpha}{1-\alpha}\ge\frac{1}{3} - \log\frac{B}{B-2} -2 \log\frac{\alpha}{1-\alpha}>0,\]
provided that $\alpha\in(10/19,1)$ is sufficiently close to $10/19$. This proves that $p_g \le x = \log^B|g|$ with $B=19$ for sufficiently large $|g|$.

\begin{rmk}\label{rmk:p_g}
Since $\tilde{A}_1(g)\ge\tilde{A}_1(21)=4/5$ for $g>1$, the proof of Theorem \ref{thm:pgupper} shows that the exponent $B=19$ can be improved to $16$ if we focus merely on positive $g\in\Gg$. Besides, if we write $g=g_1m^2$ with $g_1\in\Z$ squarefree and $m\in\N$, then $\tilde{A}_1(g)=1+o(1)$ whenever $|g|$ is sufficiently large, provided that $m^2=o(|g|)$ or $g_1\not\equiv1\pmod{4}$. Consequently, our proof of Theorem \ref{thm:pgupper} yields $p_g\ll \log^{13}(2|g|)$ for these $g\in\Gg$. In particular, this inequality holds for all squarefree $g\in\Gg$.
\end{rmk}

\section{The average value of $p_g$: Proof of Corollary \ref{cor:pgaverage}}\label{sec:avgvalue}

We remind the reader that $r$ is always to be understood as representing a prime number. We let $r_1 = 2, r_2=3, r_3=5, \dots$ denote the sequence of primes in the usual increasing order.
\begin{lem}\label{lem:vaughan} For a certain constant $\varrho >1$, we have
\[ \prod_{r \le r_k} \left(1-\frac{\phi(r-1)}{r}\right) = \varrho^{-(1+o(1))k} \quad\text{as}\quad k\to\infty. \]
\end{lem}

Lemma \ref{lem:vaughan} and the prime number theorem together imply that
\begin{equation}\label{eq:lemalternative} \prod_{r \le y} \left(1-\frac{\phi(r-1)}{r}\right) = \exp(-(1+o(1))(\log{\varrho})y/\log{y}), \end{equation}
as $y\to\infty$. We make repeated use below of this form of Lemma \ref{lem:vaughan}.

\begin{proof}[Proof of Lemma \ref{lem:vaughan}]
We will prove the lemma for a constant $\varrho$ constructed in terms of the moments of $\frac{\phi(r-1)}{r-1}$.

We start by observing that $\frac{\phi(r-1)}{r} \le \frac12$ for all primes $r$. This is clear for $r=2$, while when $r$ is odd, $\frac{\phi(r-1)}{r} < \frac{\phi(r-1)}{r-1} = \prod_{p \mid r-1}\big(1-\frac{1}{p}\big) \le \frac{1}{2}$. Now for each real $\theta$ with $|\theta|\le \frac12$, and each positive integer $M$, we have $\log(1-\theta) = -\sum_{m \le M} \frac{\theta^m}{m} + O(2^{-M})$. Thus, if we define $L_k$ by the equation $\prod_{r \le r_k} (1-\frac{\phi(r-1)}{r})=\exp(-L_k)$, then 
\[ L_k = \sum_{m \le M} \frac{1}{m} \sum_{r \le r_k} \left(\frac{\phi(r-1)}{r}\right)^m + O(2^{-M} k). \]
Here $M$ is a positive integer parameter at our disposal.

Continuing, note that $(\frac{\phi(r-1)}{r})^m - (\frac{\phi(r-1)}{r-1})^m \ll_{M} \frac{1}{r}$ for all primes $r$ and all positive integers $m\le M$. Hence, for all $k\ge 2$,
\begin{align} L_k &= \sum_{m \le M} \frac{1}{m}\sum_{r \le r_k} \left(\frac{\phi(r-1)}{r-1}\right)^m + O_M\left(\sum_{r \le r_k} r^{-1} \right) + O(2^{-M}k)\notag \\
&= \sum_{m \le M} \frac{1}{m} \sum_{r \le r_k} \left(\frac{\phi(r-1)}{r-1}\right)^m + O_M(\log\log{r_k}) + O(2^{-M} k).\label{eq:earlierexpr}
\end{align}
According to Lemma 4.4 of \cite{vaughan73}, if we set
\begin{equation}\label{eq:sigmamdef} \sigma_m := \prod_{p} \left(1-\frac{p^m-(p-1)^m}{p^{m+1}-p^m}\right),\end{equation}
then
\begin{align*} \sum_{r \le r_k} \left(\frac{\phi(r-1)}{r-1}\right)^m &= \sigma_m r_k/\log{r_k} + O_m(r_k/(\log{r_k})^2).
\end{align*}
(In \cite[Lemma 4.4]{vaughan73}, the moments of $\phi(r-1)/(r-1)$ are estimated excluding $r=2$; including $r=2$ does not change the asymptotics.) 
By the prime number theorem with the de la Vallée-Poussin error bound,
\begin{align*} \sigma_m r_k/\log{r_k} + O_m(r_k/(\log{r_k})^2) &= \sigma_m \pi(r_k) + O_m(r_k/(\log{r_k})^2)  \\
&= k \sigma_m + O_m(k/\log{k}). \end{align*}Substituting into our earlier expression \eqref{eq:earlierexpr} for $L_k$ yields
\begin{equation}\label{eq:Lkalmostdone} L_k = k\sum_{m \le M} \frac{\sigma_m}{m} + O_M(k/\log{k}) + O(2^{-M} k). \end{equation}
Inspecting the product definition \eqref{eq:sigmamdef} of $\sigma_m$, we see that $0 < \sigma_m < 2^{-m}$ for each positive integer $m$. (Note that the $p=2$ term in \eqref{eq:sigmamdef} is precisely $2^{-m}$.) It follows immediately that $\sum_{m=1}^{\infty} \frac{\sigma_m}{m}$ converges to a positive number $\varrho_0$, say. Dividing \eqref{eq:Lkalmostdone} through $k$ and  sending $k$ to infinity, we find that both  $\limsup_{k\to\infty} L_k/k$ and $\liminf_{k\to\infty} L_k/k$ are within $O(2^{-M})$ of $\sum_{m \le M} \sigma_m/m$. Sending $M$ to infinity, we conclude that $\lim_{k\to\infty} L_k/k = \rho_0$. That is, \[ L_k = (1+o(1))k\rho_0, \quad\text{as}\quad k\to\infty. \] So if we define $\varrho:= \exp(\varrho_0)$, then
\[ \prod_{r \le r_k} \left(1-\frac{\phi(r-1)}{r}\right) = \exp(-L_k) = \varrho^{-(1+o(1))k},  \]
as desired.
\end{proof}

Put $L = \log{x}/\log\log{x}$. Let $\delta_p$ be defined as in \eqref{eq:deltapdef}, and set $M_p = \prod_{r \le p} r$. 
Then $p_g = p$ precisely when $g$ belongs to one of $\delta_p M_p$ residue classes modulo $M_p$. Since $M_p \ll 3^p$,
\[ \#\{g: |g|\le x: p_g = p\} = 2\delta_p x + O(3^p). \]
As $\#([-x,x]\setminus \Gg) \ll x^{1/2}$, it follows that 
\begin{align} \notag \sum_{\substack{g \in \Gg\\ |g|\le x \\ p_g \le L}} p_g = \sum_{p \le L}p \sum_{\substack{g\in\Gg \\ |g|\le x \\ p_g = p}}1 &= \sum_{p \le L}p \left(\Bigg(\sum_{\substack{|g|\le x \\ p_g = p}}1\Bigg) + O(x^{1/2})\right)\\
&= 2x \sum_{p \le L} p\delta_p + O\bigg(\sum_{p \le L}  p(3^p +x^{1/2})\bigg) \notag\\
&= 2x \sum_{p \le L} p\delta_p + O(x^{1/2}L^2). \label{eq:prelim}\end{align}

We now extend the sum on $p$ to infinity, using Lemma \ref{lem:vaughan} to estimate the resulting error. By \eqref{eq:lemalternative},  \[ \delta_p =\frac{\phi(p-1)}{p} \prod_{r < p} \left(1-\frac{\phi(r-1)}{r}\right)\le \prod_{r \le p}\left(1-\frac{\phi(r-1)}{r}\right)
= \exp(-(1+o(1)) (\log{\varrho})p/\log{p}),
\]
where the final estimate holds as $p\to\infty$. Consequently, if we fix $c = \frac{1}{2}\log\varrho$ (for instance), then $\delta_p \ll \exp(-cp/\log{p})$ for all primes $p$, and
\[ \sum_{p > L} p \delta_p \ll \exp\big(-\frac{c}{2}L/\log{L}\big) \ll \exp(-(\log{x})^{1+o(1)}). \]
Referring back to \eqref{eq:prelim}, we deduce that
\begin{equation}\label{eq:midpoint} \sum_{\substack{g \in \Gg\\ |g|\le x \\ p_g \le L}} p_g = 2x\sum_{p} p\delta_p + O(x \exp(-(\log{x})^{1+o(1)})). \end{equation}

Next, we bound the sum of the $p_g$ taken over $g \in \Gg$, $|g|\le x$, having $p_g > L$. If $p_g > L$, then $g$ belongs to one of $\xi M$ residue classes mod $M$, where
\[ M:= \prod_{r\le L} r, \quad\text{and} \quad \xi:=\prod_{r \le L}\left(1-\frac{\phi(r-1)}{r}\right).\]
The number of such $g$ with $|g|\le x$ is $\ll \xi(x+M) \ll \xi x$, noting that $M \le 3^L = x^{o(1)}$. By another application of \eqref{eq:lemalternative}, $\xi \le \exp(-c L/\log{L})$. (All of this is being claimed for large enough values of $x$.) Thus,
\begin{equation}\label{eq:numgbound} \#\{g: |g| \le x, p_g > L\} \ll x \exp(-cL/\log{L}), \end{equation}
so that by Theorem \ref{thm:pgupper},
\begin{align*} \sum_{\substack{g \in \Gg\\ |g|\le x \\ p_g > L}} p_g &\le (\max_{\substack{g \in \Gg \\ |g| \le x}} p_g) \#\{g: |g| \le x, p_g > L\} \\&\ll (\log{x})^{19} (x \exp(-cL/\log{L}))\\ &\ll x\exp(-(\log{x})^{1+o(1)}). \end{align*}
Putting together the contributions,
\[ \sum_{\substack{g \in \Gg\\ |g|\le x}} p_g = 2x\sum_{p} p\delta_p  + O(x \exp(-(\log{x})^{1+o(1)})). \]
Corollary \ref{cor:pgaverage} follows because the function in \eqref{eq:pgaverage} is $\sum_{p} p\delta_p + O(\exp(-(\log{x})^{1+o(1)}))$.

\section{An unconditional tamed average: Proof of Theorem \ref{thm:pgavg2}}\label{sec:tamedavg}

Our main tool for this proof will be Montgomery's ``arithmetic large sieve'' inequality introduced in Section \ref{sec:explicitp_g}. Using Montgomery's sieve, Vaughan showed \cite[eq. (1.3)]{vaughan73} that for every pair of integers $M,N$ with $N>0$, we have $p_g \le N^{1/2}$ for all $g\in [M+1,M+N]$ apart from $O(N^{1/2} (\log{N})^{1-\alpha})$ exceptions, where $\alpha$ is an explicit positive constant (see \cite[eq. (1.4)]{vaughan73} for its precise definition).
Earlier Gallagher \cite{gallagher67} had shown such a result with $1$ in place of $1-\alpha$. The next proposition implies that $N^{1/2}$ can be replaced by a large power of $\log{N}$, if one is willing to slightly inflate the exponent $1/2$ on $N$ in the size of the exceptional set.

\begin{prop}\label{prop:vaughanvariant} Let $M, N \in \Z$ with $N>100$. Let $Y$ be a real number satisfying
\[ \log^2{N} \le Y \le \exp\left(\log{N} \frac{\log\log\log{N}}{\log\log{N}}\right). \]
The count of integers $g$ in $[M+1,M+N]$ with $p_g > Y$ does not exceed
\[ N^{1/2} \exp\left(O\left(\log{N} \frac{\log\log\log{N}}{\log\log{N}}\right)\right) \cdot \exp(u\log{u}),\]
where $u:= \frac{1}{2}\frac{\log{N}}{\log{Y}}$. Here the $O$-constant is absolute.
\end{prop}

Note that if $Y = \log^K{N}$ for a fixed $K \ge 1$, then the upper bound in the conclusion of Proposition \ref{prop:vaughanvariant} assumes the form $N^{\frac12(1+1/K)+o(1)}$, as $N\to\infty$. 

\begin{proof}[Proof of Proposition \ref{prop:vaughanvariant}] We may assume $N$ is sufficiently large. We apply the arithmetic large sieve from Section \ref{sec:explicitp_g} with $Q = N^{1/2}$, taking $\nu(p) = \phi(p-1)$ for $p \le Y$, and $\nu(p) = 0$ for $Y < p \le Q$. It suffices to show that with these choices of parameters, the denominator 
\begin{equation}\label{eq:Jparticular} J = \sum_{\substack{n \le N^{1/2} \\ P^{+}(n) \le Y}} \mu^2(n) \prod_{p\mid n} \frac{\phi(p-1)}{p-\phi(p-1)}
\end{equation}
in the sieve bound satisfies
\begin{equation}
\label{eq:desiredlowerJ}
 J \ge N^{1/2}\exp\left(O\left(\log{N} \frac{\log\log\log{N}}{\log\log{N}}\right)\right) \cdot \exp(-u\log{u}).
\end{equation}

Let $R$ be the number of primes $p \in [\frac{1}{2}Y,Y]$ for which the smallest prime factor of $\frac{p-1}{2}$ exceeds $Y^{1/5}$. By the linear sieve and the Bombieri--Vinogradov theorem, $R \gg Y/(\log{Y})^2$. (This application of the linear sieve is `isomorphic' to the one described at the start of \cite[Chapter 8]{DH08}. Here $1/5$ may be replaced by any constant smaller than $1/4$.)
Let $p$ be one of these $R$ primes. Then $\frac{\phi(p-1)}{p-1} = \frac{1}{2} \prod_{\ell \mid p-1,~\ell>2}(1-1/\ell) > \frac{1}{2} (1-y^{-1/5})^{4} > 2/5$ (for instance). Hence, $\frac{\phi(p-1)}{p} > \frac{1}{3}$, and $\frac{\phi(p-1)}{p-\phi(p-1)} > \frac{1}{2}$. Let $u_0 = \lfloor \log(N^{1/2})/\log{Y}\rfloor$ (so that $u_0 =\lfloor u\rfloor$, with $u$ as in the statement of Proposition \ref{prop:vaughanvariant}). By considering the contribution to the right-hand side of \eqref{eq:Jparticular} from products of $u_0$ distinct primes $p$ of the above kind, we see that $J \ge 2^{-u_0} \binom{R}{u_0}$. Now $R > Y/(\log{Y})^3 > (\log{N})^{3/2} > u_0$. Since 
\[ \binom{n}{k} = \prod_{0 \le j < k} \frac{n-j}{k-j} \ge \left(\frac{n}{k}\right)^k \] for each pair of integers $n, k$ with $n\ge k >0$, we conclude that
\[ \frac{1}{2^{u_0}} \binom{R}{u_0} \ge (R/2u_0)^{u_0} \ge (R/2)^{u_0} \exp(-u\log{u}).\]
Furthermore, using again that $R >Y/(\log{Y})^3$,
\[ (R/2)^{u_0} \ge (R/2)^{u-1} \ge Y^{u-1} (2 (\log{Y})^3)^{-u} = N^{1/2} Y^{-1} (2(\log{Y})^3)^{-u}. \]
The assumed bounds on $Y$ ensure that 
\[Y^{-1} (2(\log{Y})^3)^{-u} = \exp\left(O\left(\log{N} \frac{\log\log\log{N}}{\log\log{N}}\right)\right).\] Our desired lower estimate \eqref{eq:desiredlowerJ} then follows by combining the last two displays.
\end{proof}

\begin{proof}[Proof of Theorem \ref{thm:pgavg2}] Fix $K\ge 2$ with $\eta + \frac{1}{2}(1+1/K) < 1$. We start by estimating the contribution of $g \in \Gg$, $|g|\le x$, having $p_g \le \log^{K}{(3x)}$.

Let $L=\log{x}/\log\log{x}$. We showed in \eqref{eq:midpoint} that (as $x\to\infty$)
\[ \sum_{\substack{g \in \Gg,~|g| \le x \\ p_g \le L}} p_g = 2x \sum_{p} p\delta_p + O(x \exp(-(\log{x})^{1+o(1)})).\]Furthermore (see \eqref{eq:numgbound}), the count of $g \in \Gg$, $|g| \le x$ with $p_g > L$ is $O(x \exp(-cL/\log{L}))$, where $c = \frac{1}{2}\log\varrho>0$. Hence,
\[  \sum_{\substack{g \in \Gg,~|g| \le x \\ L < p_g \le \log^{K}(3x)}} p_g \ll x \log^{K}(3x) \exp(-cL/\log{L}) \ll x \exp(-(\log{x})^{1+o(1)}).  \]
Combining the last two displays, 
\[ \sum_{\substack{g \in \Gg,~|g| \le x \\ p_g \le \log^{K}(3x)}} \min\{p_g, x^{\eta}\} =  \sum_{\substack{g \in \Gg,~|g| \le x \\ p_g \le \log^{K}(3x)}} p_g = 2x\sum_{p} p \delta_p + O(x \exp(-(\log{x})^{1+o(1)})), \]
as $x\to\infty$.

Therefore, the proof of Theorem \ref{thm:pgavg2} will be completed once it is shown that 
\[ \sum_{\substack{g \in \Gg,~|g| \le x \\ p_g > \log^{K}(3x)}} \min\{p_g, x^{\eta}\} = o(x), \]
as $x\to\infty$. For this we apply Proposition \ref{prop:vaughanvariant}. Choose $M$ and $N$ with $M+1=-\lfloor x\rfloor$ and $M+N=\lfloor x\rfloor$; then $[M+1,M+N]$ is the set of all integers $g$ with $|g| \le x$, and $N=2\lfloor x\rfloor + 1 < 3x$. Thus, if $p_g > \log^K{(3x)}$, then $p_g > \log^K{N}$. By Proposition \ref{prop:vaughanvariant}, the number of such $g$, $|g| \le x$, is at most $x^{\frac12(1+1/K)+o(1)}$. It follows that the sum appearing in the last display is bounded above by $x^{\eta} \cdot x^{\frac12(1+1/K)+o(1)}$, which is $o(x)$ by our choice of $K$. 
\end{proof}

\begin{rmk} Our proof of Theorem \ref{thm:pgavg2} does not make essential use of the restriction $g \in \Gg$: If we average $\min\{p_g, x^{\eta}\}$ over \emph{all} integers $g$ with $|g|\le x$, the same arguments show that the limit is again $\sum_{p} p\delta_p$. For this unrestricted average, $\frac{1}{2}$ is a natural boundary for $\eta$, in that even, square values of $g$ will send the average of $\min\{p_g, x^{\frac12+\epsilon}\}$ to infinity for any fixed $\epsilon > 0$. One might hope to push $\eta$ past $\frac{1}{2}$ after restoring the condition that $g \in \Gg$, but it is not clear how to work the restriction of $g$ to $\Gg$ into the proof of a result like Proposition \ref{prop:vaughanvariant}.
\end{rmk}

\section{Almost-primitive roots}\label{sec:pgstar} 

Recall from the introduction that $g$ is called an \textsf{almost-primitive root} mod $p$ when $g$ generates a subgroup of $(\Z/p\Z)^{\times}$ of index at most $2$. Define $p_g^{*}$ analogously to $p_g$ but with ``almost-primitive root'' in place of ``primitive root.'' We then expect that 
\begin{equation}\label{eq:forall} p_g^{*} < \infty ~\,\text{for \emph{every} nonzero $g \in \Z$}.\end{equation} This seems difficult to establish unconditionally, but it can be seen to follow from GRH by a  modification of Hooley's argument. As we are not aware of a reference, we include a short GRH-conditional proof of \eqref{eq:forall} at the end of this section.

Our final main result is an upper bound on the frequency of large values of $p_g^{*}$.

\begin{thm}\label{thm:pgstarlarge} For all $x \ge 2$, there are $O(\log^3{x})$ integers $g$, $|g|\le x$, with $p_g^{*} > \log^4{x}$.
\end{thm}

Most of this section will be devoted to the proof of Theorem \ref{thm:pgstarlarge}, but we start with a few words about the application of this theorem to the average of $p_g^{\ast}$.
Put $F(p) = 1_{p>2} \phi(\frac{p-1}{2}) + \phi(p-1)$, so that $F(p)$ is the number of almost primitive roots mod $p$. Let
\[ \delta_p^{*}= \frac{F(p)}{p} \prod_{r < p} \left(1-\frac{F(r)}{r}\right). \]
Reasoning as in the introduction, we expect $p_g^{*}$ to have mean value $\sum_{p} p \delta_p^{*}$. Under GRH this could be proved analogously to our Corollary \ref{cor:pgaverage}. Using Theorem \ref{thm:pgstarlarge}, we obtain (unconditionally) that for each positive $\epsilon \in (0,1)$, the average of $\min\{p_g^{*}, x^{1-\epsilon}\}$ tends to $\sum_{p} p \delta_p^{*}$. For this,  follow the argument for Theorem \ref{thm:pgavg2} but plug in Theorem \ref{thm:pgstarlarge} in place of Proposition \ref{prop:vaughanvariant}. 

We turn now to the proof of Theorem \ref{thm:pgstarlarge}. This requires a new ingredient, Gallagher's ``larger sieve'' (see \cite{gallagher71} or \cite[\S 9.7]{FI10}). 

\begin{GLS} Let $N \in \N$, and let $\Dd$ be a finite set of prime powers. Suppose that all but $\omegab(d)$ residue classes mod $d$ are removed for each $d \in \Dd$. Then among any $N$ consecutive integers, the number remaining unsieved does not exceed
\begin{equation} \left({\sum_{d \in \Dd}\Lambda(d) - \log{N}}\right)\Bigg/\left({\sum_{d \in \Dd}\frac{\Lambda(d)}{\omegab(d)} - \log{N}}\right),\label{eq:GLSbound} \end{equation}
as long as the denominator is positive.\end{GLS}

We call $\theta \in (0,1)$ \textsf{admissible} if, for all large enough values of $Y$, we have
\[ \#\left\{p \le Y: P^{-}\left(\frac{p-1}{2}\right) > Y^{\theta}\right\} \gg \frac{Y}{\log^2{Y}}. \]
(The implied constant here is allowed to depend on $\theta$.) As remarked in the proof of Proposition \ref{prop:vaughanvariant}, the Bombieri--Vinogradov theorem in conjunction with the linear sieve implies that any $\theta < \frac14$ is admissible. It is known that there are admissible values of $\theta >\frac14$; for instance, \cite[Theorem 25.11]{FI10} shows that $\theta = \frac{3}{11}$ is admissible. 

\begin{proof}[Proof of Theorem \ref{thm:pgstarlarge}] We prove a somewhat more general result. Fix an admissible $\theta \in (0,1)$. Let $x$ be a large real number, and define
\begin{equation}\label{eq:ydef} y =((\log{x}) (\log\log{x})^2)^{1/\theta}. \end{equation}
We show that
\begin{equation}\label{eq:moregeneral} \#\{g: |g|\le x, p_g^{\ast} >y\} \ll y^{1-\theta}. \end{equation}
Theorem \ref{thm:pgstarlarge} follows from \eqref{eq:moregeneral} upon choosing an admissible $\theta > \frac14$.

We sieve the $N:=2\lfloor x\rfloor + 1$ integers in the interval $[-x,x]$. Let
\[ \mathcal{D} = \left\{\text{primes $p$}: 3<p \le y, P^{-}\left(\frac{p-1}{2}\right) > y^{\theta}\right\}. \]
Since $\theta$ is admissible,
\[ \#\Dd \gg \frac{y}{\log^2{y}}. \]
For each $p \in \Dd$, we remove every residue class \emph{except} $0\bmod{p}$ and the classes corresponding to integers whose multiplicative order mod $p$ does not exceed
\[ z:= y^{1-\theta}. \]
Then, in the notation of the larger sieve,
\begin{equation}\label{eq:omegabar}
\omegab(p) = 1 + \sum_{\substack{f \mid p-1 \\ f\le z}} \phi(f). 
\end{equation}

Suppose the integer $g$, $|g|\le x$, is removed in the sieve. In this case, there is a prime $p\in \Dd$ not dividing $g$ for which the order of $g$ mod $p$, which we will call $e$, exceeds $z$. Then
\[ \frac{p-1}{e} < \frac{y}{e} < \frac{y}{z} = y^{\theta}. \]
Since every odd prime divisor of $p-1$ exceeds $y^{\theta}$, the ratio $\frac{p-1}{e}$ cannot be divisible by any odd prime. Thus, $\frac{p-1}{e} = 2^j$ for a nonnegative integer $j$. Since $2^j\mid p-1$ and $p\equiv 3\pmod{4}$, either $j=0$ or $j=1$. That is, $e = \frac{p-1}{2}$ or $p-1$. Hence, $g$ is an almost-primitive root mod $p$. In particular, $p_g^{\ast} \le p \le y$.

Therefore, the number of $g$, $|g|\le x$, with $p_g^{\ast} > y$ is bounded above by the count of unsieved integers, which can be approached with the larger sieve. The arguments below draw inspiration from Gallagher's proof of \cite[Theorem 2]{gallagher71}.

By the Cauchy--Schwarz inequality,
\begin{equation}\label{eq:CSineq} \Bigg(\sum_{p \in \Dd} \frac{\log{p}}{\omegab(p)}\Bigg)  \Bigg(\sum_{p \in \Dd} \omegab(p)\log{p} \Bigg) \ge \Bigg(\sum_{p \in \Dd} \log{p}\Bigg)^2 \gg ((\log y)\#\Dd)^2 \gg \frac{y^2}{\log^2{y}}. \end{equation}
(We use here that $\log{p} \gg \log{y}$ for each $p \in \Dd$, which follows from $P^{-}(\frac{p-1}{2}) > y^{\theta}$.) 
On the other hand, referring back to \eqref{eq:omegabar},
\begin{align*} \sum_{p \in \Dd} \omegab(p)\log{p} &\le \sum_{p \in \Dd} \log{p} + \sum_{f \le z} \phi(f) \sum_{\substack{p \in \Dd \\ p\equiv 1\pmod*{f}}}\log{p} \\
&\ll (\log{y}) \#\Dd + \log{y} \sum_{f \le z}\phi(f) \#\{p \in \Dd: p \equiv 1\pmod*{f}\}.
\end{align*}
Brun's sieve implies that $\#\Dd \ll y/\log^2{y}$. Brun's sieve also handles the counts appearing in the sum on $f$: If $p \in \Dd$, $p\equiv 1\pmod{f}$, and $p> y^{\theta}$, then $t:=\frac{p-1}{f}< y/f$, and both $tf+1$ and $t$ have no odd prime factors up to $y^{\theta}$. Brun's sieve shows that the number of such $t$ is 
\[ \ll \frac{y}{f} \prod_{2 < r\le y^{\theta}} \left(1-\frac{1 + 1_{r\nmid f}}{r}\right) \\
\ll \frac{y}{f \log^2{y}} \prod_{r\mid f}\left(1-\frac{1}{r}\right)^{-1} = \frac{y}{\phi(f) \log^2{y}}.\]
Since there are trivially at most $y^{\theta}/f$ primes up to $y^{\theta}$ in the residue class $1\bmod{f}$, 
\[ \#\{p \in \Dd: p \equiv 1\pmod*{f}\} \ll \frac{y}{\phi(f) \log^2{y}}, \]
and 
\[ \log{y} \sum_{f \le z}\phi(f)\#\{p \in \Dd: p \equiv 1\pmod*{f}\} \ll \frac{yz}{\log{y}}. \]
We conclude that 
\[ \sum_{p \in \Dd} \omegab(p)\log{p} \ll \frac{yz}{\log{y}},\]
and hence by \eqref{eq:CSineq},
\[ \sum_{p \in \Dd} \frac{\log{p}}{\omegab(p)} \gg \frac{y^2/\log^2{y}}{yz/\log{y}} = \frac{y^{\theta}}{\log{y}}. \]
Recalling our definition \eqref{eq:ydef} of $y$, we have that
\[ \frac{y^{\theta}}{\log{y}} \gg (\log{x})(\log\log{x}), \]
which is of larger order than $\log{N}$. Hence, the denominator in \eqref{eq:GLSbound} is $\gg y^{\theta}/\log{y}$. The numerator in \eqref{eq:GLSbound} is bounded above by $\sum_{p \in \Dd} \log{p} \le (\log{y})\#\Dd \ll y/\log{y}$. Therefore, the number of unsieved $g$, $|g|\le x$, is 
\[ \ll \frac{y/\log{y}}{y^{\theta}/\log{y}} = y^{1-\theta}. \]
This completes the proof of \eqref{eq:moregeneral}.
\end{proof}

\begin{rmk} It seems likely that every $\theta \in (0,1)$ is admissible. If so, \eqref{eq:moregeneral} implies that the exponents $3$ and $4$ in Theorem \ref{thm:pgstarlarge} can be brought arbitrarily close to $0$ and $1$, respectively. \end{rmk}

\begin{proof}[Proof of \eqref{eq:forall}, assuming GRH]
Fix $g\in \Z$, $g\ne 0$. If $g \in \{\pm 1\}$, then $p_g^{\ast} = 2$. So we may assume that $|g| > 1$. As before, we let $h$ denote the largest positive integer for which $g \in (\Q^{\times})^h$. Since $g$ is fixed, we will allow implied constants below to depend on $g$ (and hence also on $h$, as fixing $g$ fixes $h$).

To prove $p_g^{\ast}$ exists, it is enough to show there is some prime $p\equiv 3\pmod{4}$, $p\nmid g$, with the property that $p$ passes the $\ell$-test for every odd prime $\ell$. Indeed, if $p$ is any such prime and $t$ is the order of $g$ mod $p$, then $\frac{p-1}{t}$ is a divisor of $p-1$ not divisible by any odd prime $\ell$. Hence, $\frac{p-1}{t}$ is a power of $2$. As $p\equiv 3\pmod{4}$, we must have $\frac{p-1}{t} = 1$ or $2$, so that $t=\frac{p-1}{2}$ or $t=p-1$. Therefore, $g$ is an almost primitive root mod $p$. (We encountered a similar argument in the proof of Theorem \ref{thm:pgstarlarge}.)

Let $x$ be large, and let $\mathcal{P} = \{\text{primes $3 < p \le x$}: P^{-}(\frac{p-1}{2}) > x^{1/5}\}$. Then $\#\mathcal{P} \gg x/(\log{x})^2$. We will show that as $x\to\infty$, all but $o(x/(\log{x})^2)$ primes $p\in \mathcal{P}$ have the desired property. In particular, there is at least one such $p$.

Clearly (once $x$ is large), each $p\in\mathcal{P}$ belongs to the residue class $3$ mod $4$. Suppose now that $p \in \mathcal{P}$ but that $p$ fails the $\ell$-test for an odd prime $\ell$. As $\ell\mid p-1$, we have $x^{1/5} < \ell \le x$.

Suppose to start with that $\ell \in (x^{1/5}, x^{1/2}/(\log{x})^3]$. The number of $p\in \mathcal{P}$ failing the $\ell$-test is certainly no more than the total number of primes $p \le x$ failing the $\ell$-test, which can be bounded by Lemma \ref{eq:failsdtest}. Indeed, using that $[\Q(\zeta_\ell,\sqrt[\ell]{g}):\Q] = (\ell-1) \frac{\ell}{(\ell,h)} \gg \ell^2$, we get from Lemma \ref{eq:failsdtest} that there are $O(\frac{x}{\ell^2\log{x}} + x^{1/2}\log{x})$ such $p$. Summing on $\ell$, the number of $p \in \Pp$ failing the $\ell$-test for some $\ell \in (x^{1/5}, x^{1/2}/(\log{x})^3]$ is $O(x/(\log{x})^3) = o(x/(\log{x})^2)$.

Suppose next that $p\in \mathcal{P}$ fails the $\ell$-test for a prime $\ell \in (x^{1/2}/(\log{x})^3, x^{1/2} (\log{x})^3]$. Write $p= 1+2\ell m$. Then $m \le x/2\ell$ and $m$ avoids the residue classes $0\bmod{r}$ and $-(2\ell)^{-1} \bmod{r}$ for each odd prime $r\le x^{1/5}$; furthermore, these two residue classes are distinct for each $r$. Noting that $x^{1/5} \le x/2\ell$, Brun's sieve bounds the number of possibilities for $m$ given $\ell$ (and hence for $p$ given $\ell$) as
\[ \ll \frac{x}{2\ell} \prod_{2<r \le x^{1/5}} \left(1-\frac{2}{r}\right) \ll \frac{x}{\ell (\log{x})^2}. \]
Summing on $\ell$ with Mertens' theorem, we conclude that the number of $p \in \Pp$ failing the $\ell$-test for some $\ell \in (x^{1/2}/(\log{x})^3, x^{1/2} (\log{x})^3]$ is $O(x \log\log{x}/(\log{x})^3) = o(x/(\log{x})^2)$.

Finally, we suppose that $p\in \mathcal{P}$ fails the $\ell$-test for an $\ell > x^{1/2}(\log{x})^3$. Then the order of $g$ mod $p$ is at most $x^{1/2}/(\log{x})^3$. Hence, $p \mid g^m-1$ for a natural number $m < x^{1/2}/(\log{x})^3$. For each $m\in \N$, the number of prime divisors of $g^m-1$ is $O(m)$. Summing on $m < x^{1/2}/(\log{x})^3$, we see that the number of $p$ arising in this way is $O(x/(\log{x})^6)$, which is certainly $o(x/(\log{x})^2)$. This completes the proof.
\end{proof}

\section*{Acknowledgments}
We thank Igor Shparlinski for a stimulating email message that prompted these investigations. We also found useful a reply left by Andrew Granville on MathOverflow \cite{MO}. The referee caught many typos in our initial submission, and we are grateful to them for their care and attention. This paper was written while P.P. was supported by NSF award DMS-2001581.

\bibliographystyle{amsplain}
\bibliography{ArtinPR}
\end{document}